\documentclass[preprint,3p,times,twocolumn]{elsarticle}
\usepackage{epstopdf}
\usepackage{graphicx}
\usepackage{amssymb}
\usepackage{tabularx}
\usepackage{subcaption} 
\usepackage{amsmath}
\usepackage{pifont}
\usepackage{supertabular}
\usepackage{eurosym}
\usepackage{xcolor} 
\usepackage{mathtools}
\usepackage{multirow}
\usepackage{appendix}
\usepackage{booktabs}
\usepackage{array}
\newcommand{\coo}{\ensuremath{\mathrm{CO_2}}}

\usepackage{amsthm}

\newtheorem{definition}{Definition}
\usepackage[ruled]{algorithm2e}
\usepackage{algpseudocode}
\usepackage{float}
\usepackage{centernot}

\usepackage[draft=False]{hyperref} 
\hypersetup{
     colorlinks=true,
     linkcolor=blue,
     filecolor=blue,
     citecolor=black,      
     urlcolor=cyan,
     }

\begin{document}

\begin{frontmatter}
\author[ULIEGE]{Antoine Dubois\corref{ADUBOIS}}
\ead{antoine.dubois@uliege.be}
\author[RTE]{Jonathan Dumas}
\author[UCLOUVAIN]{Paolo Thiran}
\author[UCLOUVAIN]{Gauthier Limpens}
\author[ULIEGE,PARIS]{Damien Ernst}
\address[ULIEGE]{Department of Computer Science and Electrical Engineering, Liege University, Liege, Belgium}
\address[RTE]{Research and Development, RTE, Paris, France}
\address[UCLOUVAIN]{Institute of Mechanics, Catholic University of Louvain, Louvain-la-Neuve, Belgium}
\address[PARIS]{LTCI, Telecom Paris, Institut Polytechnique de Paris, Paris, France}
\cortext[ADUBOIS]{Corresponding author}

\title{
Multi-objective near-optimal necessary conditions\\ for multi-sectoral planning
}

\begin{abstract}
This paper extends the concepts of \textit{epsilon-optimal spaces} and \textit{necessary conditions} for near-optimality from single-objective to multi-objective optimisation.
These notions are first presented for single-objective optimisation, and the mathematical formulation is adapted to address the multi-objective framework.
Afterwards, we illustrate the newly developed methodology by conducting multi-sectoral planning of the Belgian energy system with an open-source model called EnergyScope TD.
The cost and energy invested in the system are used as objectives.
Optimal and efficient solutions for these two objectives are computed and analysed.
These results are then used to obtain necessary conditions corresponding to the minimum amount of energy from different sets of resources, including endogenous and exogenous resources.
This case study highlights the high dependence of Belgium on imported energy while demonstrating that no individual resource is essential on its own.
\end{abstract}
\begin{keyword}
 Multi-objective optimisation; near-optimality; necessary conditions; energy system modelling; multi-sectoral planning
\end{keyword}
\end{frontmatter}

\section{Introduction} \label{sec:introduction}
Energy system planning determines the appropriate mix of energy sources and technologies to satisfy a community's or region's future energy demand.
The goal of this process is to inform decision-makers to allow them to plan an efficient and sustainable transformation of energy systems.
Energy system optimisation models (ESOMs) are commonly used to perform energy system planning \citep{decarolis2017formalizing}.
These models rely on optimisation techniques to predict how an energy system should evolve.
However, ESOMS are often used in ways that limit the quality of their insights and, thus, their usefulness to decision-makers.
Indeed, these insights are often derived from a unique cost-optimal solution.
While cost is a crucial indicator of the affordability and viability of an energy system, focusing solely on this objective might negatively impact other important factors, such as environmental sustainability and social equity.
Moreover, these insights might not satisfy stakeholders with diverging interests.

\subsection{The cost as leading indicator: limits and solutions}
\begin{figure*}[!ht]
    \centering
    \includegraphics[width=\textwidth]{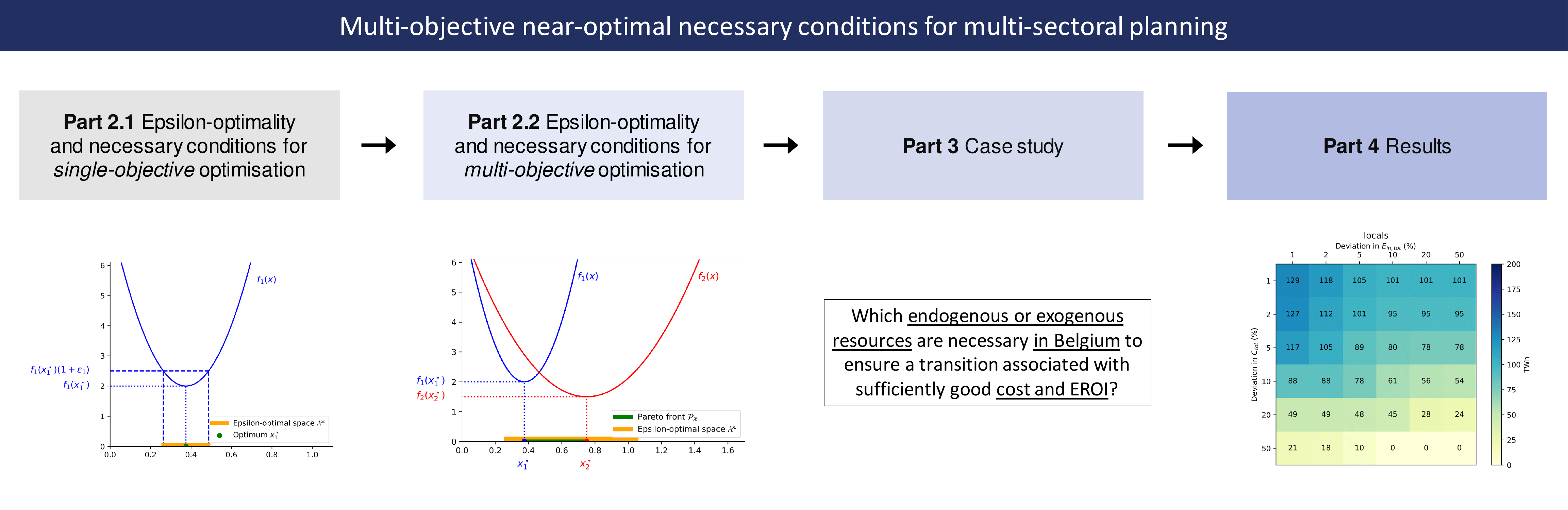}
    \caption{Graphical abstract showing the structure of the paper.
    The figures are miniatures of figures located further in the document.}
    \label{fig:graphical-abstract}
\end{figure*}
ESOMs determine the energy system configurations that minimise or maximise a specified objective.
Most studies choose the cost as the objective, and the best configuration is the most cost-effective \citep{decarolis2017formalizing}.
This choice is historical, as explained by \citet{pfenninger2014energy}.
Indeed, the first ESOMs (from the MARKAL/TIMES \citep{fishbone1981markal} and MESSAGE \citep{schrattenholzer1981energy} models) were initially designed for cost minimisation.
The study of \citet{yue2018review} highlights that by default, ESOMs ignore non-economic factors entering into energy investment decisions and how politics, social norms, and culture shape public policies.
This claim is also supported by \citet{pfenninger2014energy}, who specifies that energy system models focus heavily on economic and technical aspects.
This focus is inadequate for energy system planning as this problem involves multiple stakeholders with different policy objectives, for whom cost-optimal solutions might not be satisfying.
Moreover, several studies have demonstrated that ignoring non-economic factors increases the uncertainty of the models \citep{pfenninger2014energy, yue2018review}.
\citet{fazlollahi2012methods} also states that, due to uncertainty in some parameters, it is insufficient for energy system sizing to rivet on a unique mono-objective optimal solution.
Finally, \citet{trutnevyte2016does} shows how cost-optimal scenarios do not adequately represent real-world problems.
However, there exist methods for going beyond cost and considering non-economic factors.

\subsubsection{Scenario analysis}

The first approach to incorporate non-economic factors is scenario analysis.
Scenario analysis involves optimising the same model over multiple scenarios with different values for some parameters. 
Differences between scenarios can result from uncertainties over technological or economic parameters - e.g. future cost of technology. 
However, they can also stem from political (e.g. nuclear decommissioning) or social considerations (e.g.  limitation of onshore wind turbines or transmission lines development).
Using scenarios that differ through those considerations allows for studying the effects of non-economic factors.
For instance, the study by \citet{fujino2008back} compares a fast-growth, technology-oriented scenario to a slow-growth, nature-oriented one.
However, using scenarios to include non-economic factors is not perfect. 
Indeed, \citet{hughes2010methodological} reviewed studies using scenarios for low greenhouse gases (GHG) emission targets, and they conclude that those modelling studies still have the weakness of simplifying social and political dynamics.

\subsubsection{Multi-objective optimisation}

The second approach to include non-economic factors is multi-objective optimisation.
This approach allows for optimising several objectives simultaneously, highlighting the \textit{trade-offs} that can be obtained.

More formally, while different methods exist to apply multi-objective optimisation (e.g. weighted-sum approach, integer cut constraints, $\epsilon$-constraint method, evolutionary algorithm), they exhibit the common goal of obtaining solutions from a Pareto optimal set, also called the \textit{Pareto front}.
This Pareto front is composed of \textit{efficient} solutions, i.e. solutions that are at least better than any other solutions in one objective.
Thus, it is composed of the set of optimal trade-offs between the studied objectives, i.e. any solution that is not part of the Pareto front is worse in all objectives than at least one solution in the Pareto front.

Using multi-objective optimisation, the cost can still be optimised while considering other indicators.
For instance, \citet{becerra2008multi} conducted a study of a Texas power generation system analysing the trade-offs between economic and exergetic costs, i.e. the cumulative exergy - \textit{entropy-free energy} - consumption.
They demonstrated how these trade-offs provide insights to the decision-makers by not focusing exclusively on economic cost.
Other objectives, such as water consumption, grid dependence on imports or energy system safety, are compared to cost by \citet{fonseca2021multi, fonseca2021sustainability}.
They show how much the assessed criteria impact the design and operation of distributed energy systems.
A final example of an alternative objective that is often combined with the cost is the amount of carbon emissions \citep{jing2018multi}.

\subsubsection{Near-optimal spaces analysis} \label{sec:intro-near-optimal-spaces}

The last methodology that allows taking social and political factors into account is the study of near-optimal spaces, also called sub-optimal or epsilon-optimal spaces.
The idea is to analyse solutions close to the optimal solution to understand how the use of resources and technologies varies when allowing a slight deviation in the objective function.
This paradigm goes further than multi-objective optimisation, as mentioned by \citet{decarolis2011using}.
It allows incorporating unmodelled objectives, typical of social factors, as they are unknown or difficult to model.
Indeed, the near-optimal region might contain solutions that are worse in terms of the main objective - e.g. the cost of the system - but better in terms of unmodelled objectives such as risk or social acceptance.
This concept was introduced in the 1980s by \citet{brill1982modeling}.
The authors proposed the first method for exploring those spaces: the Hop-Skip-Jump method.
This algorithm was coined as part of a broader exploration methodology that the authors refer to as Modelling to Generate Alternatives (MGA).
This methodology was brought back recently and applied to energy system modelling by \citet{decarolis2011using} and \citet{decarolis2016modelling}.
They led to a renewed interest in such methods.
Authors such as \citet{PRICE2017356} developed new exploration algorithms while \citet{li2017investment} combined MGA with Monte-Carlo exploration to minimise parametric uncertainty.

There are several ways of extracting insights from near-optimal spaces.
Most researchers exploring near-optimal spaces focus on computing numerous near-optimal solutions from which they derive insights \cite{PRICE2017356, li2017investment, pedersen2021modeling}.
An alternative approach is to use methods to obtain such insights directly without needing to compute many alternative solutions \cite{NEUMANN2021106690}. 
The authors of \citet{DUBOIS2022108343} took this approach by introducing the concept of \textit{necessary conditions} for near-optimality, i.e. conditions that are true for every solution in the near-optimal space.
For instance, this can provide insight into the required capacity in a given technology to retain a certain level of system cost-effectiveness.
More specifically, \citet{DUBOIS2022108343} showed how, for instance, at least 200 GW of new offshore wind need to be installed Europe-wise to not deviate by more than 10\% from the cost optimum.

\subsection{Research gaps, scientific contributions and organisation}
The exploration of near-optimal spaces has been used in mono-objective optimisation problems but not, according to the best knowledge of the authors, in multi-objective optimisation problems.
However, these methods could also be valuable in multi-objective optimisation setups.
Indeed, while modelling and integrating more objectives, multi-objective optimisation still leaves aside some unmodeled objectives.
Analysing solutions in the near-optimal space of multi-objective optimisation problems is a method to address this issue. 
This paper thus aims to fill this gap by:
\begin{enumerate}
    \item extending the concepts related to near-optimal spaces to multi-objective optimisation;
    \item computing necessary conditions in a multi-objective context to highlight the range of insights that can be derived from them.
\end{enumerate}
The first point is addressed in Section~\ref{sec:problem-statement} by first introducing the mathematical concepts of near-optimality and necessary conditions in a single-objective framework (see Section~\ref{sec:methodology-single}) and then extending them to multi-objective optimisation (see Section~\ref{sec:methodology-multi}).
Section~\ref{sec:case-study} then translates those concepts to a real case study: the multi-sectoral expansion of the Belgian energy system.
The results of this case study, including necessary conditions representing the necessary amount of different energy resources, are presented in Section~\ref{sec:results} before highlighting the contributions of this paper in Section~\ref{sec:conclusion}.
We can already highlight one of those contributions: the open-source release of the code \citep{dubois2023code} and the data \citep{dubois2023data} used to achieve this study.
The graphical representation of the organisation of this paper is depicted in Figure~\ref{fig:graphical-abstract}.

\section{Problem statement} \label{sec:problem-statement}
The first part of this section introduces the concepts of epsilon-optimal space and necessary conditions for single-objective optimisation \citep{DUBOIS2022108343}.
The second part extends these concepts to multi-objective optimisation by:
\begin{enumerate}
    \item generalising the optimisation problem to multiple objectives,
    \item presenting generic notions related to multi-objective optimisation, including the image of the feasible set, efficient solutions, and the Pareto front, and
    \item explaining the extension of the concepts of epsilon-optimality and necessary conditions to multi-objective optimisation.
\end{enumerate}

\subsection{Epsilon-optimality and necessary conditions for single-objective optimisation} \label{sec:methodology-single}

\begin{figure}
    \centering
    \includegraphics[width=\linewidth]{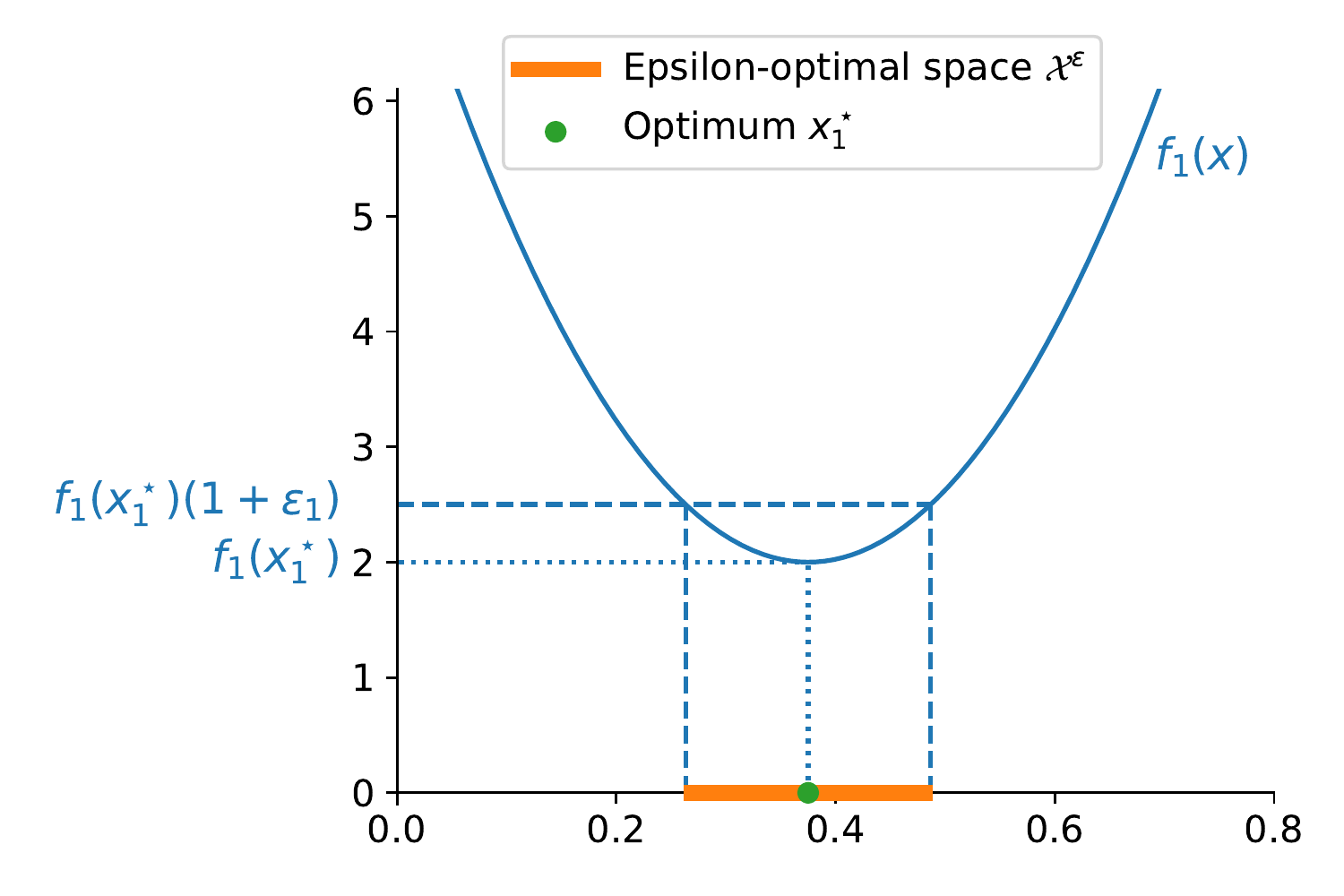}
    \caption{
    Graphical representation of an $\epsilon$-optimal space of a mono-objective optimisation problem in $\mathcal{X} = \mathbb{R}_+$. 
    The function $f_1$ that is minimised is shown in blue.
    Its minimum is located at $x_1^{\star}$.
    Using this value and its corresponding objective value $f_1(x_1^*)$ allows to easily determine an $\epsilon$-optimal space $\mathcal{X}^\epsilon$ with $\epsilon = \epsilon_1$.
    The values of the different parameters and functions used in this example are described in \ref{sec:appendix-example}.
    }
    \label{fig:epsilon-optimal-space-uni}
\end{figure}

Let $\mathcal{X}$ be a feasible space and $f:\mathcal{X} \rightarrow \mathbb{R}_{+}$ an objective function in the positive reals. 
The single-objective optimisation problem is
\begin{align} \label{equ:original-problem}
    \min_{x\in\mathcal{X}} f(x)\quad.
\end{align}
Let $x^\star$ denote an optimal solution to this problem that is: $x^\star \in \arg \min_{x\in\mathcal{X}} f(x)$.
\begin{definition}
    An \textbf{$\epsilon$-optimal space}, with $\epsilon \geq 0$, is defined as follows
    \begin{align} \label{equ:epsilon-space}
        \mathcal{X}^\epsilon = \bigg\{x \in \mathcal{X} \mid f(x) \leq (1+\epsilon) f(x^\star) \bigg\}\quad.
    \end{align}
\end{definition}
It is the set of the feasible solutions $x \in \mathcal{X}$ with objective value $f(x)$ no greater than $(1+\epsilon)f(x^\star)$. 
The deviation from the optimal objective value is measured via $\epsilon$, called the \textit{suboptimality coefficient}.
Figure~\ref{fig:epsilon-optimal-space-uni} illustrates those concepts.
A note must be made on the specific case $f(x^\star) = 0$.
In this case, $\mathcal{X}^\epsilon$ resumes to $\arg \min_{x\in\mathcal{X}} f(x)$, making the analysis of near-optimal spaces trivial.

\subsubsection{Necessary conditions}

The concepts of \textit{condition}, \textit{necessary condition} and \textit{non-implied necessary condition} introduced in this section allow determining features which are common to all solutions in a given $\epsilon$-optimal space.
\begin{definition}
    A \textbf{condition} is a function $\phi: \mathcal{X} \rightarrow \{0, 1\}$. 
    A set of conditions is denoted $\Phi$.
\end{definition}
\noindent We illustrate this definition and the following using an example. 
Let $\mathcal{X} = \mathbb{R}$, $\Phi$ can be the set of conditions of the form $\phi_c(x) \coloneqq x \ge c$ with $x \in \mathcal{X}$ (thus $x\in\mathbb{R})$ and $c \in \mathbb{R}$.

\begin{definition} \label{def:nec-cond-1}
    A \textbf{necessary condition} for $\epsilon$-optimality is a condition which is true for any solutions in $\mathcal{X}^\epsilon$.
    
    For a given feasible space $\mathcal{X}$, set of conditions $\Phi$ and suboptimality coefficient $\epsilon$, $\phi \in \Phi$ is a necessary condition if $$\forall x \in \mathcal{X}^\epsilon: \phi(x)=1\quad .$$
    The set of all necessary conditions for $\epsilon$-optimality in $\Phi$ is denoted $\Phi^{\mathcal{X}^\epsilon}$.
\end{definition}
\noindent In our example, let $\mathcal{X}^\epsilon = [0, 1]$, then $\phi_0(x) \coloneqq x \ge 0$ is respected by all $x \in \mathcal{X}^\epsilon$, making $\phi_0$ a necessary condition.
Moreover, the set of all necessary conditions is $\Phi^{\mathcal{X}^\epsilon} = \{\phi_c \mid c\leq 0 \}$.
Indeed, it is straightforward to show that $\phi_c(x) \coloneqq x \ge c$ is true over $\mathcal{X}^\epsilon = [0, 1]$ for any $c\leq 0$.

As shown in \citet{DUBOIS2022108343}, necessary conditions can provide insights into features common to many near-optimal solutions.
However, depending on how conditions are defined, their study also claims the number of necessary conditions can be infinite, which is counterproductive in providing insights.
This is, for instance, the case in our previous example.
Indeed, the set $\Phi^{\mathcal{X}^\epsilon} = \{\phi_c \mid c\leq 0 \}$ contains an infinite number of necessary conditions.
Nevertheless, the only condition of interest is $\phi_0$.
Indeed, knowing that $x \geq 0$ is true over $[0, 1]$ implies that $x \geq c$ is true when $c<0$.
Thus, knowing that $\phi_0$ is a necessary condition implies that any $\phi_c$ with $c<0$ is a necessary condition.  
On the opposite, it is not possible to imply that $\phi_0$ is a necessary condition from the knowledge of other necessary conditions in the set $\Phi^{\mathcal{X}^\epsilon} = \{\phi_c \mid c\leq 0 \}$.
This defines $\phi_0$ as a \textit{non-implied} necessary condition.
Knowing this, we can provide the first definition of non-implied necessary conditions.
\begin{definition}
A condition $\phi \in \Phi$ is a \textbf{non-implied necessary condition} if it cannot be implied to be a necessary condition from the sole knowledge of other necessary conditions.
The set of non-implied necessary conditions is denoted $\overline{\Phi}^{\mathcal{X}^\epsilon}$.
\end{definition}

\noindent This definition can be made more formal.
In the following, we formalise it mathematically, starting with the notion of \textit{implication}.
\begin{definition}
An \textbf{implication function} $\psi(\phi_1\mid \phi_2) \in \{0, 1\}$ is a function that indicates whether condition $\phi_2$ implies condition $\phi_1$.
When $\psi(\phi_1\mid \phi_2) = 1$, then $\forall x \in \mathcal{X},\; \phi_2(x) = 1 \implies \phi_1(x) = 1$.
When  $\psi(\phi_1\mid \phi_2) = 0$, then $\exists x \in \mathcal{X},\; \phi_2(x) = 1 \centernot\implies \phi_1(x) = 1$.
\end{definition}
\noindent Using this function, a non-implied necessary condition can be defined in the following way.
\begin{definition}
    A \textbf{non-implied necessary condition} is a necessary condition $\phi \in \Phi^{\mathcal{X}^\epsilon}$ that is not implied by any other necessary condition.
    It is a necessary condition $\phi \in \Phi^{\mathcal{X}^\epsilon}$ such that $\forall \phi' \in \Phi^{\mathcal{X}^\epsilon}\setminus\{\phi\}: \psi(\phi\mid \phi') = 0$.
\end{definition}
\noindent We can refine the definitions of necessary and non-implied necessary conditions by defining the space over which a condition is true.
\begin{definition}
    The space $\mathcal{I}_\phi$ is the subset of $\mathcal{X}$ where a condition $\phi$ is true, that is:
    \begin{align}
        \mathcal{I}_\phi = \bigg\{x\in\mathcal{X} \mid \phi(x)=1 \bigg\}\quad.
    \end{align}
\end{definition}
\noindent Necessary conditions can then be defined in the following way.
\begin{definition}
    A condition $\phi$ is a \textbf{necessary condition} for $\epsilon$-optimality if $\mathcal{X}_\epsilon \subseteq \mathcal{I}_\phi$.
\end{definition}

\begin{figure}[!ht]
    \centering
    \includegraphics[width=1.5in]{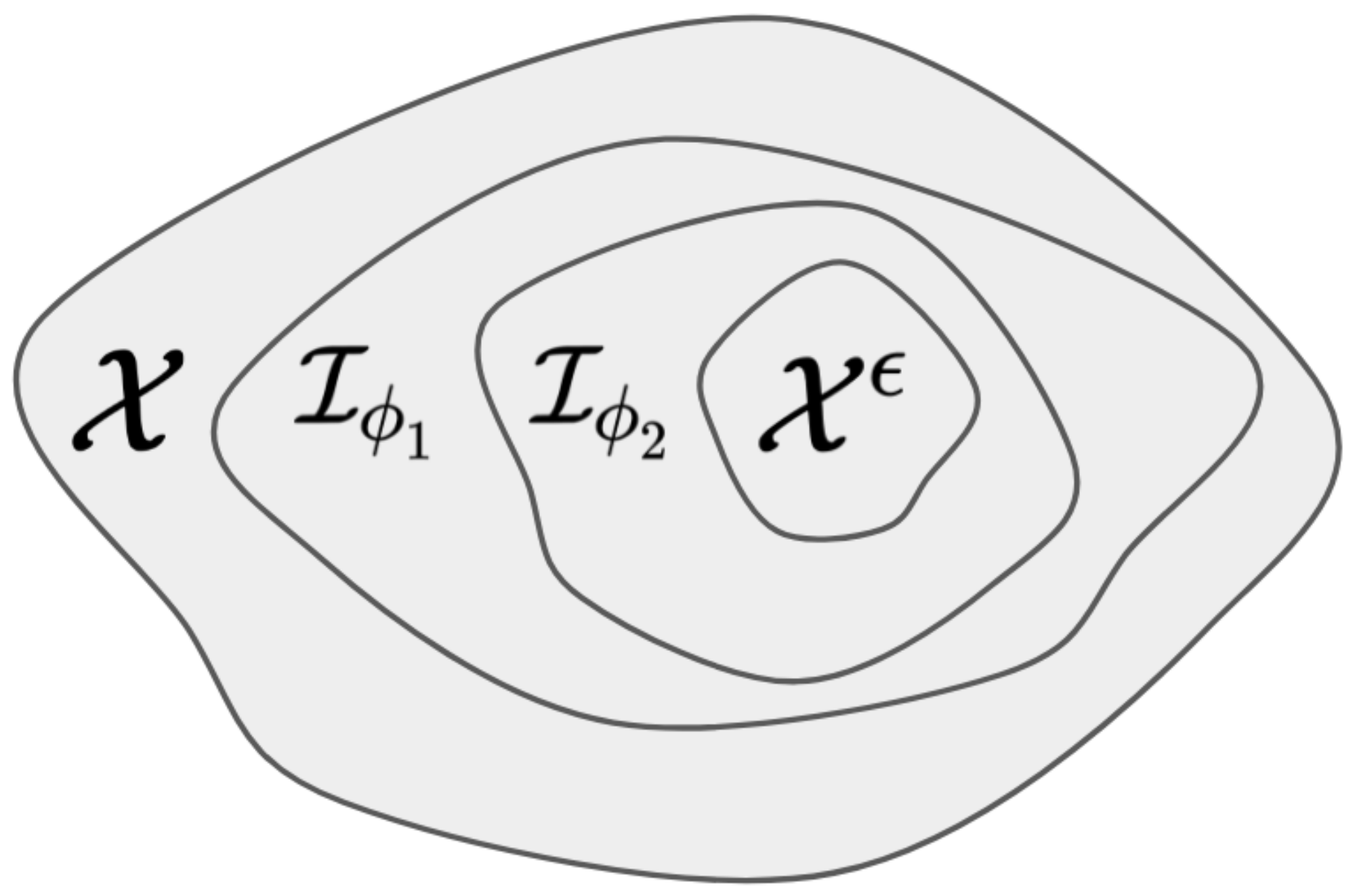}%
    \caption{Graphical illustration of implication using spaces over which conditions are true, adapted from \citet{DUBOIS2022108343}.
    The four spaces that are represented are the feasible space $\mathcal{X}$, the $\epsilon$-optimal space $\mathcal{X}^\epsilon$, and the spaces $\mathcal{I}_{\phi_1}$ and $\mathcal{I}_{\phi_2}$ where conditions $\phi_1$ and $\phi_2$ are respectively true.
    Both these conditions are necessary as $\mathcal{X}^\epsilon \subset \mathcal{I}_{\phi_1}$ and $\mathcal{X}^\epsilon \subset \mathcal{I}_{\phi_2}$.
    Moreover, $\phi_2$ implies $\phi_1$ as $\mathcal{I}_{\phi_2} \subset \mathcal{I}_{\phi_1}$.
    }
    \label{fig:true-spaces}
\end{figure}

\noindent Indeed, this implies that $\forall x \in \mathcal{X}_\epsilon,\; \phi(x) = 1$, which corresponds to definition \ref{def:nec-cond-1}.
Using this same notion, we can give an implementation of the implication function.
\begin{definition}
    Let $\phi_1$ and $\phi_2$ be conditions with $\mathcal{I}_{\phi_1}$ and $\mathcal{I}_{\phi_2}$ the spaces over which they are respectively true, then the \textbf{implication function} $\psi(\phi_1 \mid \phi_2)$ is defined as:
    \begin{align}
        \psi(\phi_1 \mid \phi_2) = \mathcal{I}_{\phi_2} \subseteq \mathcal{I}_{\phi_1}\quad.    
    \end{align}
\end{definition}
\noindent This formulation fits our definition of an implication function because if $\psi(\phi_1 \mid \phi_2) = 1$, then it means $\mathcal{I}_{\phi_2} \subseteq \mathcal{I}_{\phi_1}$, which in turns implies that $\phi_1(x) = 1$ for any $x \in \mathcal{X}$ for which $\phi_2(x) = 1$.
Similarly, if $\psi(\phi_1 \mid \phi_2) = 0$, it means that $\mathcal{I}_{\phi_2} \not\subseteq \mathcal{I}_{\phi_1}$, which means $\exists x \in \mathcal{X}$ such that $\phi_1(x) = 0$ when $\phi_2(x) = 1$.
This definition of implication leads to a new definition of non-implied necessary conditions.
\begin{definition}
    A \textbf{non-implied necessary condition} is a necessary condition $\phi \in \Phi^{\mathcal{X}^\epsilon}$ that is true over a space which does not include any of the spaces over which other necessary conditions are true.
    It is a necessary condition $\phi \in \Phi^{\mathcal{X}^\epsilon}$ such that $\forall \phi' \in \Phi^{\mathcal{X}^\epsilon}\setminus\{\phi\}: \mathcal{I}_{\phi'} \not\subseteq \mathcal{I}_{\phi}$.
\end{definition}
\noindent Figure~\ref{fig:true-spaces} illustrates these concepts, where $\phi_2$ implies $\phi_1$ as $\mathcal{I}_{\phi_2} \subset \mathcal{I}_{\phi_1}$. 
They are both necessary conditions because they are true over $\mathcal{X}^\epsilon$.
Finally, if there are no other conditions in the set of conditions $\Phi = \{\phi_1, \phi_2\}$, then $\phi_2$ is a non-implied necessary condition as no other necessary condition implies it.

\subsubsection{Computation of a non-implied necessary condition} \label{sec:example-non-implied-nec-cond}

This section presents an example to provide a practical sense of these concepts. 
It demonstrates how to compute a non-implied necessary condition from a set of conditions taking the form of constrained sums of variables. 
In the case studies described in Section \ref{sec:case-study}, this type of condition is used to study the minimum amount of energy that can be driven from different sources.

Let $\mathcal{X} \subset \mathbb{R}^n$ be a feasible space, $f:\mathcal{X} \rightarrow \mathbb{R}_{+}$ an objective function to minimise over this space, and $\Phi_{\mathbf{d}}$ a set of conditions defined as follows:
\begin{align} 
    \Phi_{\mathbf{d}} = \bigg\{&\phi_{\mathbf{d}}^ c(\mathbf{x}) = \mathbf{d}^T\mathbf{x} \geq c \bigg\}\quad,
\end{align}
where $\mathbf{x} \in \mathcal{X}$, $\mathbf{d} \in \{0, 1\}^n$ and $c \in \mathbb{R}$.
The conditions are constrained sums of variables $\mathbf{d}^T\mathbf{x}=\sum_{i=1}^n d_ix_i$.
In this particular case, \citet{DUBOIS2022108343} have proven that $\phi_\mathbf{d}^{c^\star}= \mathbf{d}^T\mathbf{x} \geq c^*$ with $c^* = \min_{\mathbf{x}\in\mathcal{X}^\epsilon} \mathbf{d}^T\mathbf{x}$ is the only non-implied necessary condition that can be derived from $\Phi_{\mathbf{d}}$.
The value $c^*$ represents the minimum value that $\mathbf{d}^T\mathbf{x}$ can take over the set $\mathcal{X}^\epsilon$, that is when allowing a deviation of $\epsilon$ from the optimal value $f(\mathbf{x}^\star)$.
Algorithm \ref{alg:cap-1} illustrates the computation of this value in three steps.
\begin{algorithm}
\caption{Computation of a non-implied necessary condition - Single-objective case}\label{alg:cap-1}
    \KwData{\\\quad$f$ - objective function,\\
    \quad $\mathcal{X}$ - feasible space,\\
    \quad $\epsilon$ - suboptimality coefficient,\\
    \quad $\mathbf{d}$ - binary vector defining the conditions $\mathbf{d}^T\mathbf{x}$\\}
    \KwResult{$c^*$}
    \textbf{Steps:}
    \begin{enumerate}
        \item Solve $\min_{\mathbf{x}\in \mathcal{X}} f(\mathbf{x})$ to obtain $\mathbf{x}^\star$.
        \item Build $\mathcal{X}^\epsilon$ by adding to the original problem the constraint $f(\mathbf{x}) \le (1+\epsilon)f(\mathbf{x}^\star)$.
        \item Solve $c^* = \min_{\mathbf{x}\in\mathcal{X}^\epsilon}\mathbf{d}^T\mathbf{x}$.
    \end{enumerate}
\end{algorithm}

\subsection{Epsilon-optimality and necessary conditions for multi-objective optimisation} \label{sec:methodology-multi}

This section extends the concepts presented previously to multi-objective optimisation while introducing notions specific to this type of optimisation problem.

Let $\mathbf{f} := (f_1, \cdots, f_k, \cdots, f_n)$ be a vector of $n$ objective functions such that $\forall k\; f_k:\mathcal{X} \rightarrow \mathbb{R}_{+}$. 
We seek to minimise these functions over the feasible set $\mathcal{X}$, which, using the notation of \citet{ehrgott2005multicriteria}, we note: 
\begin{align} \label{equ:multi-criteria-initial-problem}
    ``\min_{x \in \mathcal{X}}"\quad \mathbf{f}\quad.
\end{align}

\noindent Let $\mathcal{Y}$ be the image of the feasible set in the objective space:
\begin{align}
    \mathcal{Y} = \mathbf{f}(\mathcal{X}) = \{y \in \mathbb{R}^n\mid y = \mathbf{f}(x) \text{ for some } x \in \mathcal{X}\}\quad.
\end{align}
This space is the image of $\mathcal{X}$ under the objective functions $\mathbf{f}$, and $\mathbf{f}(x) := (f_1(x), \cdots, f_k(x), \cdots, f_n(x))$. 
Therefore, $\mathcal{Y} \in \mathbb{R}^n_{+}$ and each of its components $y_k$ are defined by $y_k = f_k(x)$ for some $x \in \mathcal{X}$.
\begin{figure*}[!ht]
    \centering
    \subfloat[Full representation of $\mathcal{X}^{\boldsymbol{\epsilon}}$]
    {\includegraphics[width=3in]{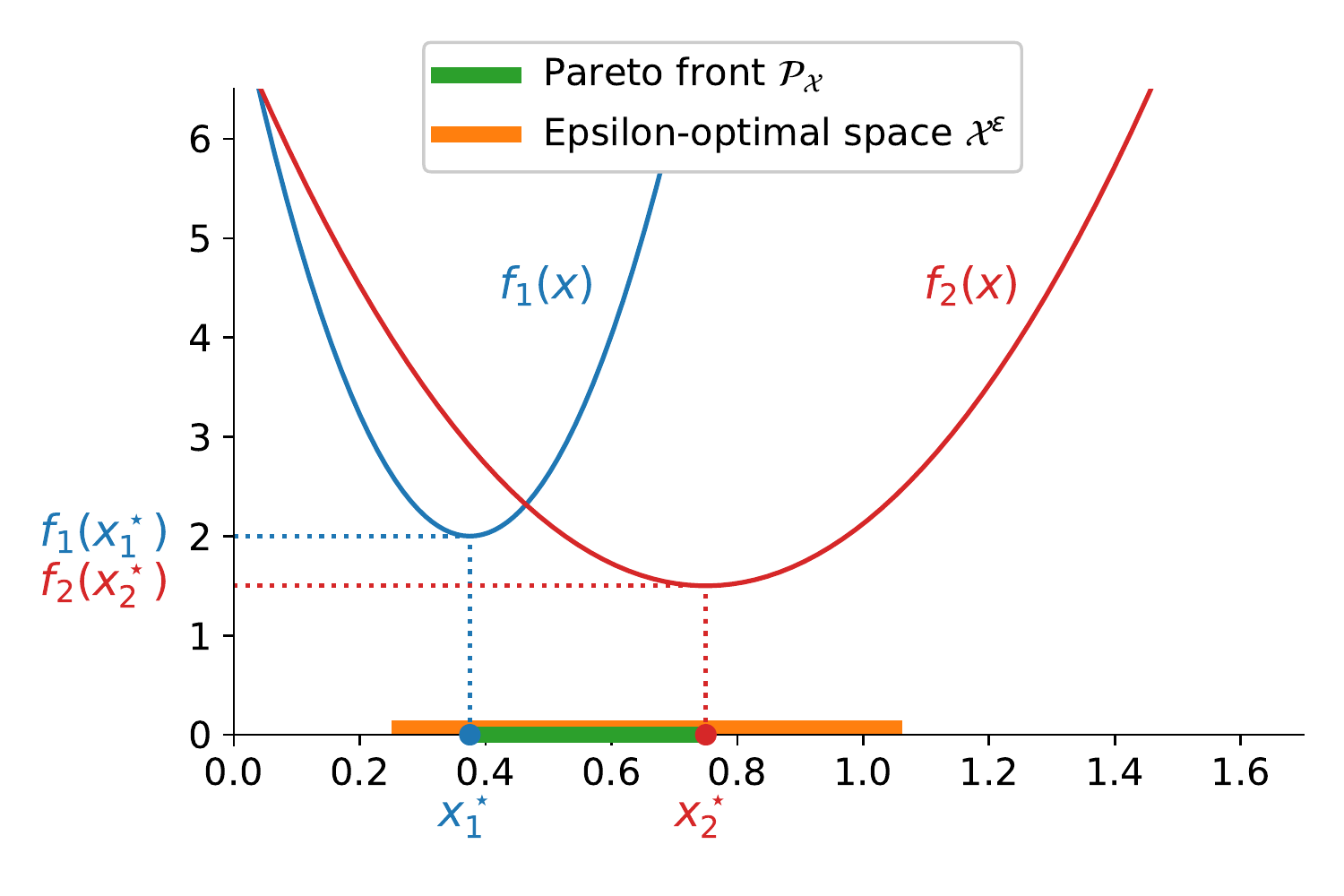}%
    \label{fig:epsilon-optimal-space-2d-full}}
    \hfil
    \subfloat[Computation of a subset of $\mathcal{X}^{\boldsymbol{\epsilon}}$ from a single point $\hat{x}$ of the Pareto front $\mathcal{P}_\mathcal{X}$]
    {\includegraphics[width=3in]{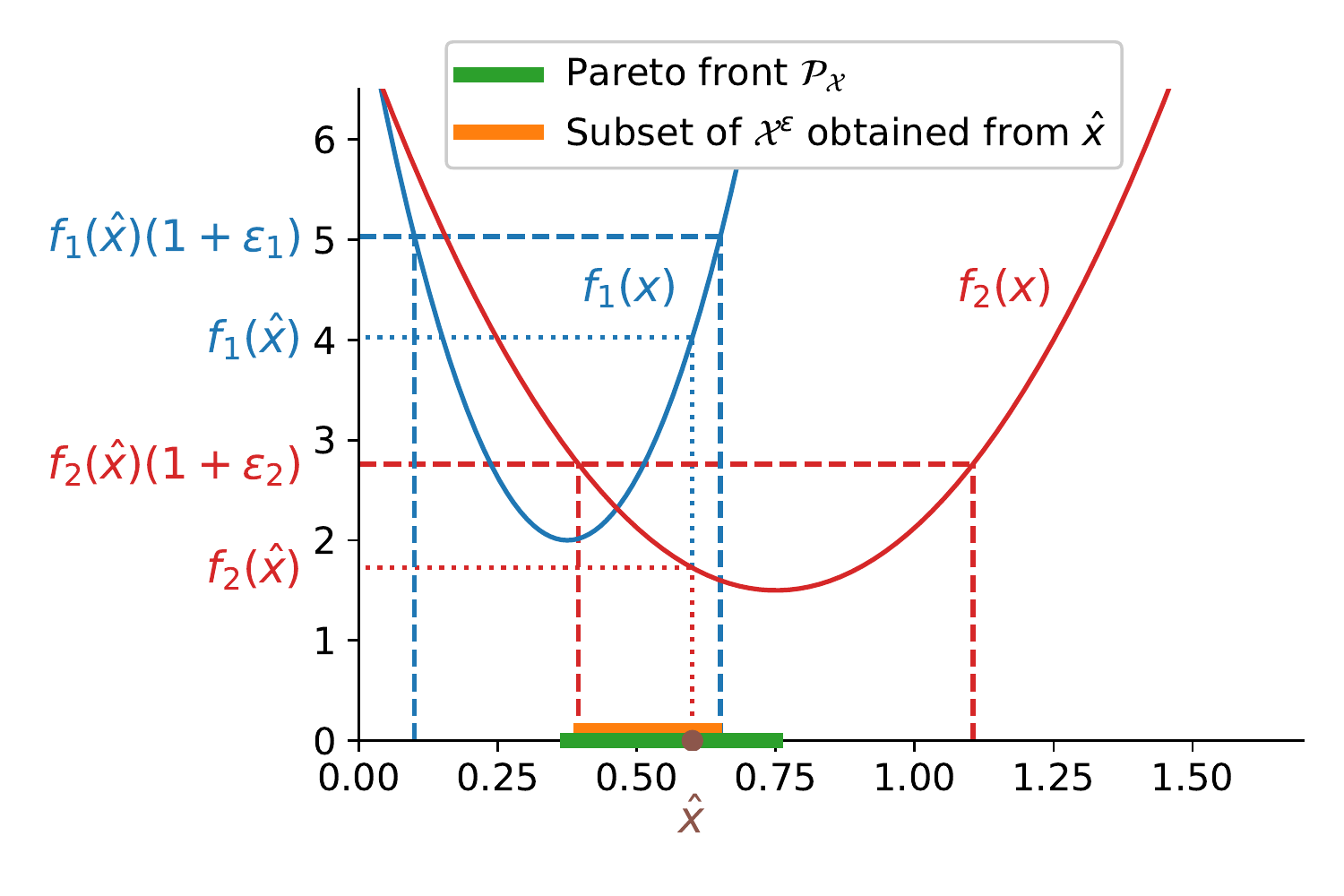}%
    \label{fig:epsilon-optimal-space-2d-part}}
    \caption{
        Graphical representation of an $\boldsymbol{\epsilon}$-optimal space of a multi-objective optimisation problem in $\mathcal{X} = \mathbb{R}_+$.
        The two functions to be minimised $f_1$ and $f_2$ are represented in blue and red, respectively, and their respective minimums are $x_1^{\star}$ and $x_2^{\star}$.
        The Pareto front $\mathcal{P}_\mathcal{X}$ containing all efficient solutions is represented in green.
        Figure \ref{fig:epsilon-optimal-space-2d-full} shows in orange the full $\boldsymbol{\epsilon}$-optimal space $\mathcal{X}^{\boldsymbol{\epsilon}}$ for a suboptimality coefficient vector $\boldsymbol{\epsilon} = (\epsilon_1, \epsilon_2)$. 
        As shown in Equation~\eqref{equ:epsilon-space-multi-criteria-def2-2}, this space is the union of sub-spaces that can be computed from efficient solutions.
        Figure \ref{fig:epsilon-optimal-space-2d-part} shows how one of these subspaces, corresponding to the efficient solution $\hat{x}$, can be computed.
        From the value $\hat{x}$, the corresponding objective values $f_1(\hat{x})$ and $f_2(\hat{x})$ are obtained.
        This allows to determine all the solutions in $\mathcal{X}$ whose objective value is smaller than $f_k(\hat{x})(1+\epsilon_k)$ for $k \in {1, 2}$.
        The values of the different parameters and functions used in this example are described in \ref{sec:appendix-example}.
    }
    \label{fig:epsilon-optimal-space-2d}
\end{figure*}

\subsubsection{Efficient solutions and Pareto front}
One way to highlight compromises between the objectives ($f_1, \cdots, f_k, \cdots, f_n$) is to compute efficient or Pareto optimal solutions.
Following the definition of \citet{ehrgott2005multicriteria}:
\begin{definition}
    A feasible solution $x \in \mathcal{X}$ is called \textbf{efficient} when there is no other $\hat{x} \in \mathcal{X}$ such that $\forall k\; f_k(\hat{x}) \leq f_k(x)$ and $f_i(\hat{x}) < f_i(x)$ for some $i$, that is, no other $\hat{x} \in \mathcal{X}$ has a smaller or equal value in all objectives $(f_1, \cdots f_k, \cdots, f_n)$ than $x$.     
\end{definition}
According to \citet{ehrgott2005multicriteria}, multiple denominations exist for the set of efficient points. 
This paper uses the term Pareto front to indiscriminately name the set of efficient points or their image in the objective space.
\begin{definition} \label{def:pareto-front}
    A \textbf{Pareto front} $\mathcal{P}_{\mathcal{X}}$ is the set
    \begin{align}
        \mathcal{P}_{\mathcal{X}} = \bigg\{x\in \mathcal{X}\mid & \not\exists \hat{x} \in \mathcal{X},\\\nonumber 
        & \forall k\; f_k(\hat{x}) \leq f_k(x), \exists i\; f_i(\hat{x}) < f_i(x) \bigg\}\quad.
    \end{align}
    In the objective space, a Pareto front is defined as:
    \begin{align}
        \mathcal{P}_{\mathcal{Y}} = \bigg\{y\in \mathcal{Y}\mid & \not\exists \hat{y} \in \mathcal{Y},\\\nonumber 
        &\forall k\; \hat{y}_k \leq y_k, \exists i\; \hat{y}_i < y_i \bigg\}\quad.
\end{align}
\end{definition}

\subsubsection{Approximation of Pareto fronts} \label{sec:appr-pareto-front}

A Pareto front can be composed of an infinity of points. 
Thus, it is typical to compute a subset of the efficient solutions which compose it.
This set is named \textit{approximated Pareto front} and is denoted by $\mathcal{P}_{\mathcal{X}, m}$ (or equivalently $\mathcal{P}_{\mathcal{Y}, m}$) where $m$ is the number of points in the approximation.
\begin{definition}
    An \textbf{approximate Pareto front} $\mathcal{P}_{\mathcal{X}, m}$, with $m \in \mathbb{N}$, is a subset of $m$ efficient solutions in the Pareto front $\mathcal{P}_{\mathcal{X}}$.
\end{definition}
\noindent Several techniques exist to obtain those efficient solutions, the two most famous being the `\textit{weighted-sum approach}' and the `\textit{$\epsilon$-constraint method}' \citep{ehrgott2005multicriteria}.
The weighted-sum approach consists of solving:
\begin{equation} \label{equ:method-weighted}
    \min_{x \in \mathcal{X}} \sum_{k=1}^n \lambda_kf_k(x)\quad,
\end{equation}
with $\forall k\; \lambda_k > 0$.
%
%
The $\epsilon$-constraint method resolves in solving:
\begin{equation} \label{equ:method-epsilon-constraint}
\begin{split}
    \min_{x \in \mathcal{X} } \hspace{0.1cm} & f_j(x)\\
    \text{s.t. } & f_k(x) \leq \epsilon_k \text{ for } k = 1, \cdots, n \text{ and } k \not= j\quad,
\end{split}
\end{equation}
where $\forall k\; \epsilon_k \in \mathbb{R}$. 
%

\subsubsection{Multi-criteria epsilon-optimal spaces} \label{sec:multi-criteria-optimal-spaces}

Starting from a Pareto front $\mathcal{P}_{\mathcal{X}}$, it is possible to define an $\boldsymbol{\epsilon}$-optimal space, given a \textit{suboptimality coefficients vector} of deviations in each objective: $\boldsymbol{\epsilon} = (\epsilon_1, \cdots, \epsilon_k, \cdots, \epsilon_n) \in \mathbb{R}^n_{+}$.
This space is denoted by $\mathcal{X}^{\boldsymbol{\epsilon}}$ in the decision space and $\mathcal{Y}^{\boldsymbol{\epsilon}}$ in the objective space.

In the mono-objective setup, the $\epsilon$-optimal space is defined as the set of points $x \in \mathcal{X}$ whose objective value $f(x)$ do not deviate by more than an $\epsilon$ fraction from the optimal objective value, i.e. $f(x) \leq (1+\epsilon)f(x^\star)$.
In a multi-objective case, there is no optimum but a set of efficient points composing the Pareto front.
This leads us to define the $\boldsymbol{\epsilon}$-optimal space as follows:
\begin{definition}
    In a multi-objective optimisation problem, the $\boldsymbol{\epsilon}$\textbf{-optimal space} $\mathcal{X}^{\boldsymbol{\epsilon}}$, with $\boldsymbol{\epsilon} = (\epsilon_1, \cdots, \epsilon_k, \cdots, \epsilon_n) \in \mathbb{R}^n_{+}$, is the set of points $x$ whose objective values $f_k(x)$ for each $k$ do not deviate by more than an $\epsilon_k$ fraction from the objective value $f_k(\hat{x})$ of \textit{at least one} solution $\hat{x}$ of the Pareto front $\mathcal{P}_\mathcal{X}$.\\
    
    \noindent It is the space
    \begin{align}
        \mathcal{X}^{\boldsymbol{\epsilon}} = \bigg\{x \in \mathcal{X} \mid& \exists \hat{x} \in \mathcal{P}_{\mathcal{X}}, \nonumber\\ & \forall k\; f_k(x) \leq (1+\epsilon_k) f_k(\hat{x}) \bigg\}\quad.
        \label{equ:epsilon-space-multi-criteria-def2}
    \end{align}
    Alternatively, this space can be defined as:
    \begin{align} 
        \mathcal{X}^{\boldsymbol{\epsilon}} = \bigcup_{\substack{\hat{x} \in \mathcal{P}_{\mathcal{X}}}} \bigg\{x \in \mathcal{X} \mid& \forall k\; f_k(x) \leq (1+\epsilon_k) f_k(\hat{x}) \bigg\}\quad.
        \label{equ:epsilon-space-multi-criteria-def2-2}
    \end{align}
\end{definition}
\noindent Figure~\ref{fig:epsilon-optimal-space-2d} depicts a graphical representation of an $\boldsymbol{\epsilon}$-optimal space in a multi-objective framework and how it is built from efficient solutions.

Definition (\ref{equ:epsilon-space-multi-criteria-def2}) relies on the entire Pareto front. 
However, practically, only a subset $\mathcal{P}_{\mathcal{X}, m}$ of $m$ efficient points of the Pareto front is computed and used to obtain an approximation of the $\boldsymbol{\epsilon}$-optimal space, denoted $\mathcal{X}^{\boldsymbol{\epsilon}}_m$.
\begin{definition} \label{def:epsilon-space-multi-criteria-def2-approximate}
    An approximation $\mathcal{X}^{\boldsymbol{\epsilon}}_m$, with $m \in \mathbb{N}$, of an $\boldsymbol{\epsilon}$-optimal space $\mathcal{X}^{\boldsymbol{\epsilon}}$ is the space
    \begin{align} 
        \mathcal{X}^{\boldsymbol{\epsilon}}_m = \bigg\{x \in \mathcal{X} \mid& \exists \hat{x} \in \mathcal{P}_{\mathcal{X}, m},\nonumber\\
        & \forall k\; f_k(x) \leq (1+\epsilon_k) f_k(\hat{x}) \bigg\}\quad.
        \label{equ:epsilon-space-multi-criteria-def2-approximate}
    \end{align}
    Alternatively, this space can be defined as:
    \begin{align} 
        \mathcal{X}^{\boldsymbol{\epsilon}}_m = \bigcup_{\substack{\hat{x} \in \mathcal{P}_{\mathcal{X}, m},}} \bigg\{x \in \mathcal{X} \mid& \forall k\; f_k(x) \leq (1+\epsilon_k) f_k(\hat{x}) \bigg\}\quad.
        \label{equ:epsilon-space-multi-criteria-def2-approximate-2}
    \end{align}
\end{definition}
The alternative formulation defines $\mathcal{X}^{\boldsymbol{\epsilon}}_m$ as a union of spaces, where each space is the set of points whose objective value in each $f_k$ does not deviate by more than an $\epsilon_k$ fraction from the objective values $f_k(\hat{x})$ of one solution $\hat{x}$ in the approximated Pareto front $\mathcal{P}_{\mathcal{X}, m}$.

\begin{figure*}[!ht]
    \centering
    \subfloat[Few well-spread points]
    {\includegraphics[width=0.3\textwidth]{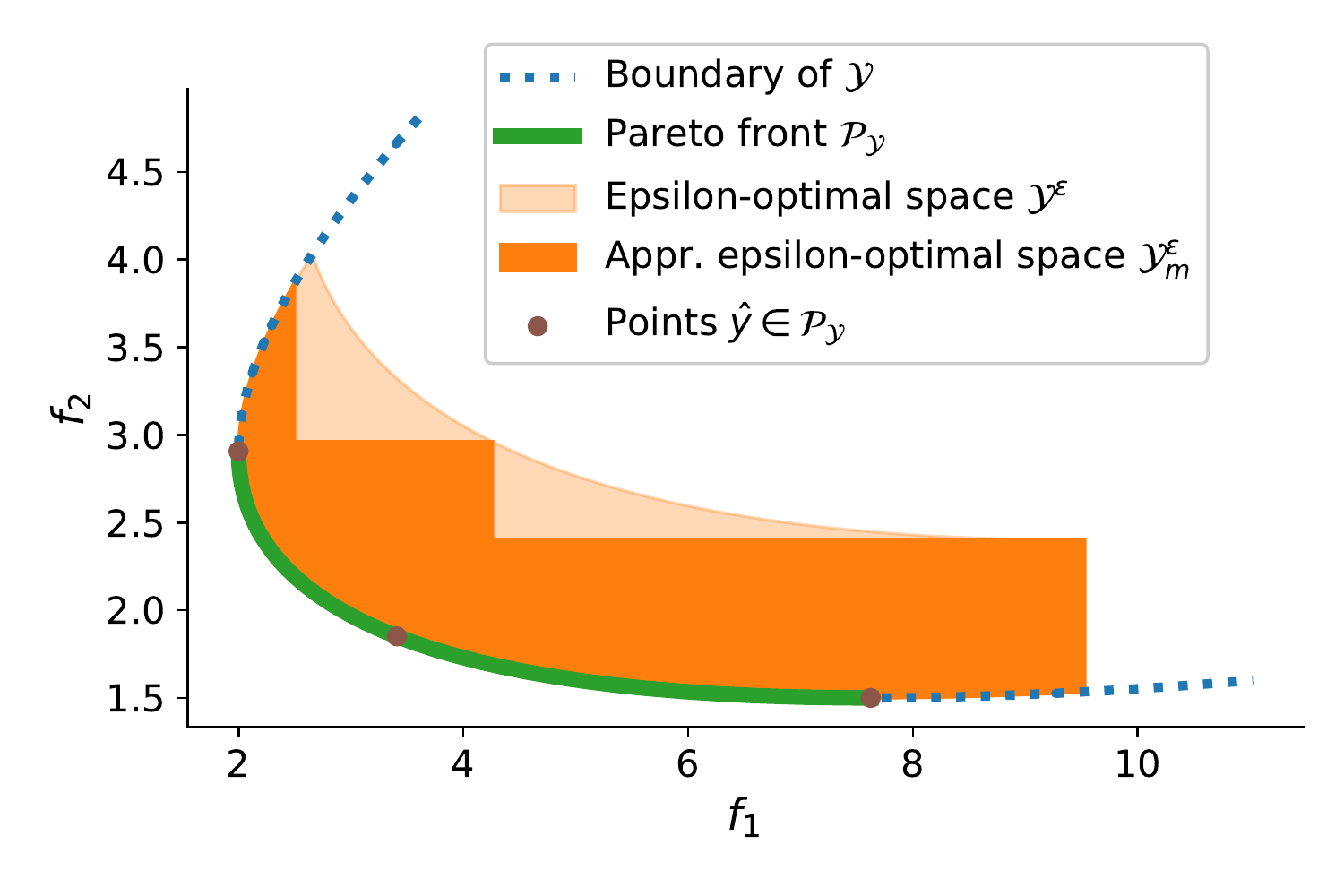}%
    \label{fig:epsilon-space-approximate-pareto-3}}
    \hfil
    \subfloat[Numerous badly-spread points]
    {\includegraphics[width=0.3\textwidth]{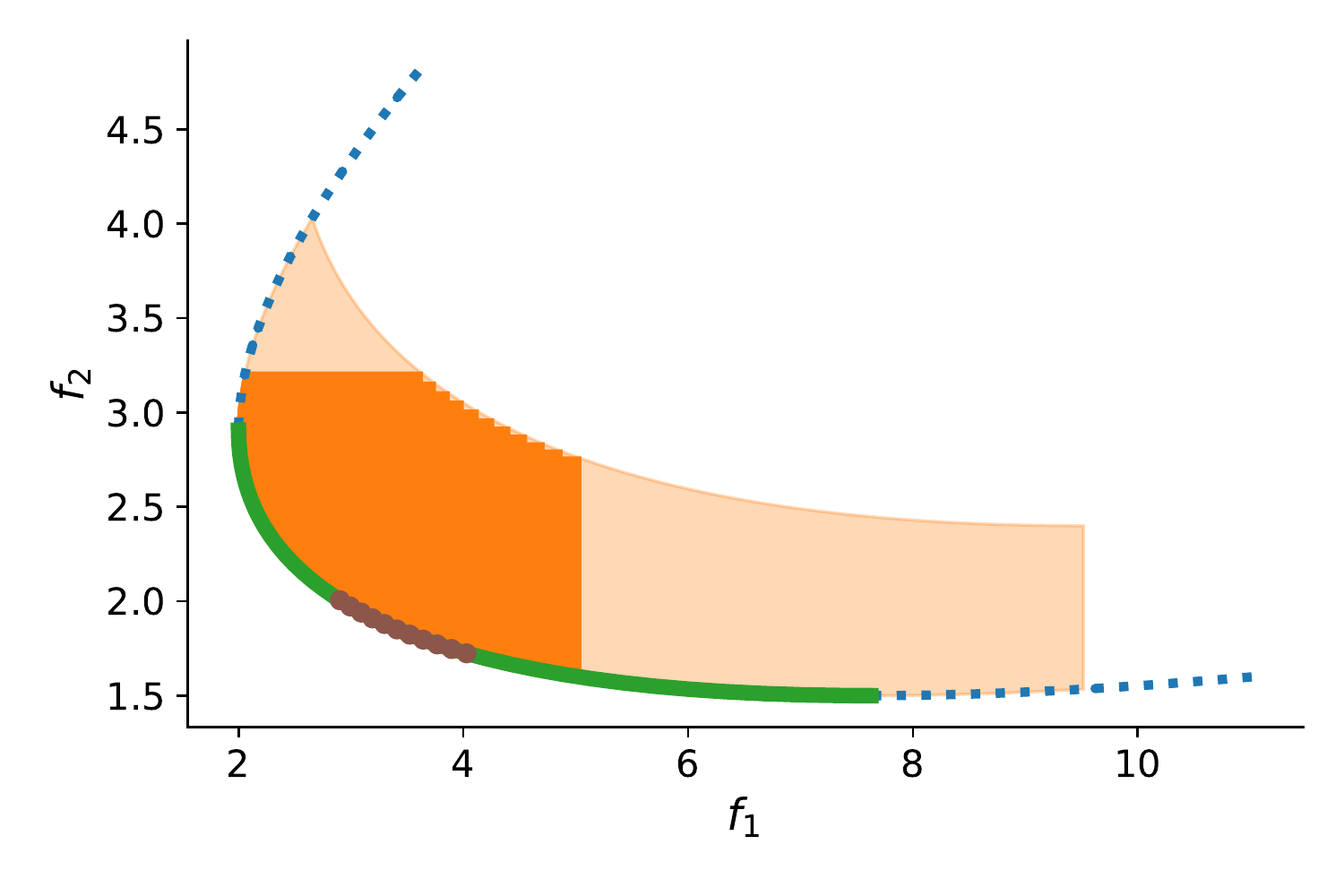}%
    \label{fig:epsilon-space-approximate-pareto-2}}
    \hfil
    \subfloat[Numerous well-spread points]
    {\includegraphics[width=0.3\textwidth]{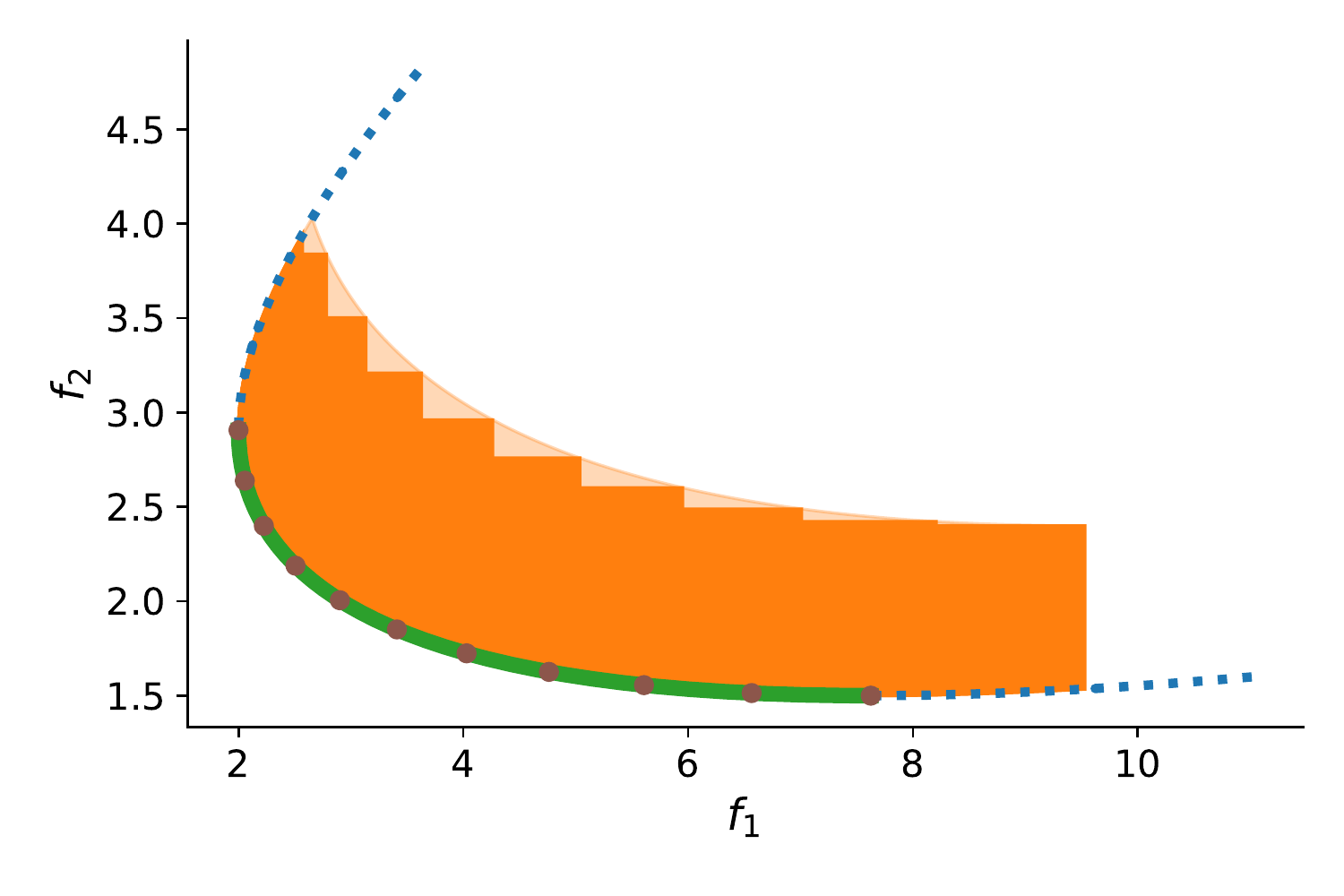}%
    \label{fig:epsilon-space-approximate-pareto-1}}
    \caption{Graphical representations in the objective space of approximations $\mathcal{Y}^{\boldsymbol{\epsilon}}_m$ of an $\boldsymbol{\epsilon}$-optimal space of a multi-objective optimisation problem based on three different approximate Pareto front $\mathcal{P}_{\mathcal{Y}, m}$.
    The axes correspond to the two functions to minimise, i.e. $f_1$ and $f_2$.
    The boundary of the image of the feasible space $\mathcal{Y}$ is represented in blue in the three cases.
    The part of this boundary corresponding to the full Pareto front $\mathcal{P}_\mathcal{Y}$ is drawn in green.
    The full $\boldsymbol{\epsilon}$-optimal spaces $\mathcal{Y}^{\boldsymbol{\epsilon}}$ corresponding to this Pareto front is coloured in light orange.
    Each graph corresponds to a different approximate Pareto front $\mathcal{P}_{\mathcal{Y}, m}$.
    These sets of points are represented in brown.
    From each of these points, part of the approximate $\boldsymbol{\epsilon}$-optimal spaces can be computed, and their union is represented in solid orange.
    The values of the different parameters and functions used in this example are described in \ref{sec:appendix-example}.
    }
    \label{fig:epsilon-space-approximate-pareto}
\end{figure*}

Figure~\ref{fig:epsilon-space-approximate-pareto} shows three examples of approximate $\boldsymbol{\epsilon}$-optimal spaces $\mathcal{X}^{\boldsymbol{\epsilon}}_m$ in the objective space (therefore noted $\mathcal{Y}^{\boldsymbol{\epsilon}}_m$) using three approximated Pareto fronts $\mathcal{P}_{\mathcal{Y}, m}$, with different numbers and spread of efficient solutions.

\subsubsection{Necessary conditions}

In the multi-objective optimisation framework, necessary conditions for $\boldsymbol{\epsilon}$-optimality can be defined in the same manner as in the one-dimensional setting, i.e. conditions which are true for any solutions in $\mathcal{X}^{\boldsymbol{\epsilon}}$. 
\begin{definition}
A \textbf{necessary condition} for $\epsilon$-optimality is a condition which is true for any solutions in $\mathcal{X}^{\boldsymbol{\epsilon}}$.
Let $\mathcal{X}$ be a feasible space, $\Phi$ a set of conditions, and $\boldsymbol{\epsilon}$ a suboptimality coefficients vector, a necessary condition $\phi \in \Phi$ is a condition such that $$\forall x \in \mathcal{X}^{\boldsymbol{\epsilon}}: \phi(x)=1\quad.$$
The set of all necessary conditions for $\boldsymbol{\epsilon}$-optimality in $\Phi$ is denoted $\Phi^{\mathcal{X}^{\boldsymbol{\epsilon}}}$.     
\end{definition}
\noindent Similarly, a non-implied necessary conditions $\boldsymbol{\epsilon}$-optimality can be defined in the following way.
\begin{definition}
    A \textbf{non-implied necessary condition} is a necessary condition $\phi \in \Phi^{\mathcal{X}^{\boldsymbol{\epsilon}}}$ that is not implied by any other necessary condition.
    It is a necessary condition such that $\forall \phi' \in \Phi^{\mathcal{X}^{\boldsymbol{\epsilon}}}\setminus\{\phi\}: \psi(\phi\mid \phi') = 0$.
\end{definition}

\subsubsection{Computation of a non-implied necessary condition} \label{sec
:computation-non-implied-necessary-conditions-multi}

The computation of a non-implied necessary condition from conditions of type $\mathbf{d}^T\mathbf{x} \geq c$ presented in Section \ref{sec:example-non-implied-nec-cond} is generalised to the multi-criteria case.
In the mono-objective case, it was sufficient to minimise the sum $\mathbf{d}^T\mathbf{x}$ over $\mathcal{X}^\epsilon$ to obtain the value $c^*$ corresponding to the non-implied necessary condition $\mathbf{d}^T\mathbf{x} \geq c^*$.
However, in a multi-objective setup, we do not have access to $\mathcal{X}^{\boldsymbol{\epsilon}}$ but to its approximation $\mathcal{X}^{\boldsymbol{\epsilon}}_m$, which is the union of several subsets, each corresponding to one point in $\mathcal{P}_{\mathcal{X}, m}$ (i.e. a subset of the Pareto front).
The minimum over this space can thus be obtained by taking the minimum of the minima of $\mathbf{d}^T\mathbf{x}$ over each of these subsets.
Even with this approach, $\mathcal{X}^{\boldsymbol{\epsilon}}_m$ being a subset of $\mathcal{X}^{\boldsymbol{\epsilon}}$, minimising $\mathbf{d}^T\mathbf{x}$ over it will only provide an upper bound $\Tilde{c}$ of the value $c^*$, i.e. $\Tilde{c} \geq c^*$.
Algorithm~\ref{alg:cap-2} shows how this value can be obtained.

\begin{algorithm}[!hb]
    \caption{Computation of a non-implied necessary condition - Multi-objective case}\label{alg:cap-2}
    \KwData{\\\quad$\mathbf{f}$ - objective functions,\\
    \quad $\mathcal{X}$ - feasible space,\\
    \quad $m$ - number of points,\\
    \quad $\boldsymbol{\epsilon}$ - vector of suboptimality coefficients,\\
    \quad $\mathbf{d}$ - binary vector defining the conditions $\mathbf{d}^T\mathbf{x}$}
    \KwResult{$\tilde{c}$}
    \textbf{Steps:}
    \begin{enumerate}
        \item Draw $m$ points $\hat{x}_1, \ldots \hat{x}_i, \ldots \hat{x}_m$ of the Pareto front using an appropriate method.
        \item For all $i \in [1, 2, \ldots, m]$, compute $c_i = \min \mathbf{d}^T\mathbf{x}$ over the space $\{x \in \mathcal{X} \mid \forall k\; f_k(x) \leq (1+\epsilon_k) f_k(\hat{x}_i)\}$.
        \item Take the minimum $\tilde{c} = \min_{i \in [1, 2, \ldots, m]} c_i$ of these values to find the appropriate condition $\phi_{\tilde{c}}$.
    \end{enumerate}
\end{algorithm}

There is no guarantee that the condition $\mathbf{d}^T\mathbf{x} \geq \tilde{c}$ is a (non-implied) necessary condition.
Indeed, it could be the case that for a solution $\mathbf{x} \in \mathcal{X}^{\boldsymbol{\epsilon}} \setminus\mathcal{X}^{\boldsymbol{\epsilon}}_m$ that $\mathbf{d}^T\mathbf{x} < \tilde{c}$.
To make the upper bound $\Tilde{c}$ as close as possible to the real minimal value $c^*$, one must reduce the size of the difference $\mathcal{X}^{\boldsymbol{\epsilon}} \setminus\mathcal{X}^{\boldsymbol{\epsilon}}_m$.
This can be done by improving the number and spread of efficient solutions in the approximated Pareto front.
As defined by \citet{alarcon2010multi}, solutions with a good spread can be seen as having good coverage of the actual Pareto front.
The three graphs of Figure~\ref{fig:epsilon-space-approximate-pareto} show visually how, by increasing the number and the spread of efficient solutions drawn from the Pareto front, the approximated $\boldsymbol{\epsilon}$-optimal space covers a more significant subset of the points of the entire $\boldsymbol{\epsilon}$-optimal space.

\section{Case study} \label{sec:case-study}

In this section, a case study will illustrate the concepts and methodology presented in the previous section.
First, the context of the case study and the question to which it tries to provide an answer are presented.
The modelling tool used to implement the methodology is then introduced, and its main features are detailed.
Finally, each element introduced in Section~\ref{sec:problem-statement} is specified to the problem at hand.

\subsection{Context} \label{sec:case-study-context}

In the European Green Deal \citep{EU-green-deal}, the European Commission raided the European Union's ambition to reduce GHG emissions to at least 55\% below 1990 levels by 2030. 
Then by 2050, Europe aims to become the world's first carbon-neutral continent.
Europe still relies massively on fossil fuels to satisfy its energy consumption ($\sim$ 75\% coming from coal, natural gas and oil according to the \citet{IEA-Data}) as well as non-energy usages (e.g. chemical feed-stocks, lubricants and asphalt for road construction \citep{doi/10.2833/001137}).
The use of these fuels is the primary source of GHG emissions.
Carbon-neutral sources of energy must thus be developed to curb emissions.
The possibilities are numerous, and one of the coming decade's main challenges will be deciding which resources to invest in.
Several criteria will motivate these choices.\\

The most common criterion for discriminating between options is cost.
Indeed, as highlighted by \citet{pfenninger2014energy} and \citet{decarolis2011using}, most studies use the cost indicator to plan the energy transition. 
This choice makes sense as the cost of investment, maintenance and operation of the energy system impacts the final consumers' energy bill.
Thus, minimising the system cost is a social imperative to allow every citizen access to affordable energy.

A lesser-known indicator, encompassing technical and social challenges, is the system's \textit{energy return on investment} (EROI).
When defined system-wise, the EROI is a ratio that measures the usable energy delivered by the system ($E_{out}$) over the amount of energy required to obtain this energy ($E_{in}$) \citep{dupont2021estimate}.
When the amount of energy required to deliver a given energy service increases, the EROI of the system decreases.
In some sense, EROI measures the ease with which energy is extracted to transform it into a form that benefits society.
There are various manners of defining $E_{in}$ and $E_{out}$, and incidentally, the EROI of a system.
These definitions depend mainly on what parts of the \textit{energy cascade} - as presented in \citet{Brockway2019} or \citet{dumas2022energy} - are considered.
This paper considers that invested energy $E_{in}$ encompasses the energy used to build the system infrastructure, `from the cradle to the grave', and to operate this system.
Following the methodology of \citet{dumas2022energy}, $E_{out}$ will correspond to the final energy consumption (FEC) of the system, as defined in the  European Commission \citep{eurostat-glossary} standard. 
FEC is the total energy, measured in TWh, consumed by end-users.
It encompasses the energy directly used by the consumer and excludes the energy used by the energy sector, e.g. deliveries and transformation.

While cost and EROI can be linked (e.g. the transport of energy resources will increase both the system cost and invested energy), they are not fully correlated and favouring one or the other can lead to different system configurations, as illustrated later in Section~\ref{sec:results}.
Both criteria can be included in the decision process by modelling them as objectives in optimisation problems.
These objectives can be optimised individually or co-optimised using multi-criteria optimisation techniques.
In this case study, we will show how, using these objectives in the methodology presented in Section~\ref{sec:problem-statement}, the following question can be addressed:
\begin{center}
    \fbox{\begin{minipage}{0.95\linewidth}
        Which resources are necessary to ensure a transition associated with sufficiently good cost and EROI?
    \end{minipage}}
\end{center}
Indeed, the answer to this question can be obtained by computing necessary conditions corresponding to the minimum amount of energy that needs to come from these resources.

This question is, however, relatively broad, and for the sake of conciseness, it needs to be specified.
On top of decision criteria, considerations such as energy independence (i.e. which has been enhanced with the Russian invasion of Ukraine) and social acceptance (i.e. the `not-in-my-backyard' phenomena) are paramount in planning the energy transition.
These considerations will impact the type of resources that will be exploited.
Indeed, the first consideration incentives a push for domestically produced energy, while the latter favours the opposite.
The first tends to minimise the amount of exogenous resources in the system, while the latter minimises the amount of energy coming from endogenous resources.
To take these elements into consideration, the previous question can be refined to:
\begin{center}
\fbox{\begin{minipage}{0.9\linewidth}
\textbf{Which endogenous or exogenous resources are necessary} to ensure a transition associated with sufficiently good cost and EROI?
\end{minipage}}
\end{center}

%
This study focuses on one of the European countries: Belgium.
Belgium made the same commitments for 2030 and 2050 as the European Union \citep{belgium-55}.
Thus, it faces the challenge of replacing its fossil-based economy with carbon-free solutions while striking the right balance between endogenous and exogenous resources.
Belgium's population density exacerbates this challenge.
In 2019, Belgium had the second-highest population density in Europe (excluding Malta) with 377 people per km$^2$, behind the Netherlands (507 people per km$^2$) \citep{eurostat-density}.
The available land for onshore energy development is thus limited, while offshore production is limited to around 8 GW of wind potential \citep{be-off-platform}.
Other domestic resources such as solar, biomass, waste or hydro also have limited potential.
This situation entails a small local energy potential compared to its demand. 
The study \citet{limpens2020belgian} evaluates that available local Belgian resources can only cover 42\% of the country's primary energy consumption.
This situation strongly impacts the type of resources Belgium must rely on.
Therefore, the question that will be addressed in this case study is:
\begin{center}
\fbox{\begin{minipage}{0.9\linewidth}
Which endogenous or exogenous resources are necessary \textbf{in Belgium} to ensure a transition associated with sufficiently good cost and EROI?
\end{minipage}}
\end{center}

\subsection{EnergyScope TD}

%
To answer this question, an appropriate ESOM is needed.
The commitments set for 2035 and 2050 cover all sectors of the economy, not just electricity production.
To achieve net zero ambitions, carbon-neutral solutions must be implemented for electricity, heat, mobility, and non-energy.
These different sectors can be modelled using an open-source whole-energy system model called EnergyScope TD (ESTD) \citep{limpens2019energyscope}.
ESTD can be categorised as an ESOM.
Using optimisation techniques, it determines the optimal mix of technologies (e.g. wind turbines, gas power plants, boilers, etc.) and resources (e.g. wind, gas, diesel, etc.) which are needed to satisfy different types of end-use demand (EUD) listed in \ref{sec:appendix-eud}.
This optimal mix depends on the user-defined objective function.
Mathematically, ESTD models the energy system as a linear programming problem.
It takes a series of parameters as input and outputs the values of investment and operational variables determined by minimising an objective while respecting a series of constraints.
The objective is a linear function; constraints are linear equalities or inequalities.
Parameters and variables can be indexed temporally.
The default temporal horizon $T$ is one year with an hourly resolution.
To reduce the computational burden of the optimisation, the horizon is clustered by selecting a number of typical days, 12 by default.
Thus, time-dependent parameters and variables are indexed by a typical day $td$ and an hour $h$.
The equivalence between the original hourly-resolution temporal horizon and the typical days is done via a time-indexed set $THTD(t)$ associating each hour $t$ of the year with a corresponding couple $(td, h) = THTD(t)$.
This set is essential to understand some of the equations in the rest of this section.

ESTD has been extensively used and validated in the Belgian case \cite{dumas2022energy, limpens2020belgian,limpens2020impact,rixhon2021role, colla2021optimal,limpens2021system}.
%
More specifically, in \citet{limpens2020belgian}, the authors studied the 2035 Belgian energy system using ESTD and built the corresponding data set.
This year is a trade-off between a long-term horizon where policies can still be implemented and a horizon short enough to define the future of society with a group of known technologies. 
To build on these resources, we will model the Belgian energy system for 2035.

To finish this section, it is important to note that while the results presented in this paper are valid for Belgium, they could easily be extended to other countries.
Indeed, ESTD has already been used to model the energy systems of other countries such as Switzerland \citep{li2020decarbonization, guevara2020machine} and Italy \citep{borasio2022deep}.
Moreover, adapting those models to implement the methodology presented in this paper only requires minor modifications, as presented in the following sections.

\subsection{Feasible space}

In the initial optimisation problem
\begin{align} \label{equ:multi-criteria-initial-problem-copy}
    ``\min_{x \in \mathcal{X}}"\quad \mathbf{f}\quad,
\end{align}
the first element to define is the feasible space $\mathcal{X}$ over which the optimisation is performed.
This study modelled the feasible space using ESTD as a linear programming problem.
Therefore, the problem to solve has the following form:
\begin{align} \label{equ:multi-criteria-initial-problem-modified}
    ``\min_{\mathbf{x}}" &\quad \mathbf{f}\\\nonumber
    s.t. &\quad A\mathbf{x} \geq \mathbf{b}\quad,
\end{align}
where $\mathbf{x}$ is the vector of variables of the problem, while $A$ and $\mathbf{b}$ are a matrix and vector of parameters, respectively.
More information on the specific variables, parameters and constraints used in ESTD can be found in \citet{limpens2019energyscope} and in the model's documentation \citep{readthedocs}.

\subsubsection{Constraint on GHG emissions}

A constraint that is of particular interest given the context of this case study is the limit on GHG emissions, i.e. 
\begin{align}
    GWP_{tot} \le 35\; \text{[Mt\coo-eq/y]}\quad.
\end{align}

In this section, we briefly describe how this constraint is defined.
%
The total yearly GHG emissions of the system are computed using a life-cycle analysis (LCA) approach. 
Thus, they include the GHG emissions along the whole life cycle, i.e. `from the cradle to the grave' of the technologies and resources considered in ESTD.
In ESTD, the global warming potential (GWP) expressed in Mt\coo-eq./year is used as an indicator to aggregate emissions of different GHG.
Then, the yearly emissions of the system, which are denoted $GWP_{tot}$, are defined as follows:
\begin{align}
    GWP_{tot} = \sum_{j\in TECH} \frac{GWP_{constr}(j)}{lifetime(j)} + \sum_{i\in RES} GWP_{op}(i)\quad,
\end{align}
where $TECH$ and $RES$ are the sets of technologies and resources modelled in ESTD.
$GWP_{constr}$ represents the GWP for the construction of a technology, while $GWP_{op}$ gives the GWP linked to the operation of a resource.
More specifically, $GWP_{constr}(j)$ is the GWP of technology $j$ over its entire lifetime allocated to one year based on the technology lifetime $lifetime(j)$.
$GWP_{op}(i)$ is the GWP related to the use of resource $i$ over one year.

The 35 Mt\coo-eq/y limit chosen in this case study comes from the following reasoning.
According to the International Energy Agency (IEA), Belgium's 1990 territorial GHG emissions were approximately 105 Mt\coo-eq \citep{IEA-Data}.
Thus, the targets of the European Green Deal imply reaching 47 Mt\coo-eq/y in 2030 and 0 Mt\coo-eq/y in 2050\footnote{Practically, the 2050 target is to be climate neutral, meaning the GHG emission can be greater than 0 but must be compensated by carbon capture.}.
By conducting a linear interpolation between these dates, the 2035 Belgian GHG emissions should reach approximately 35 Mt\coo-eq/y. 
This target is used as a hard constraint for $GWP_{tot}$ in the model: $GWP_{tot} \le 35$ [Mt\coo-eq/y].

\subsection{Objectives} \label{sec:case-study-objectives}

The second step in formalising the problem consists in choosing appropriate objectives. 
As mentioned at the start of this section, our interest lies in solutions with a sufficiently good cost and EROI.
This choice implies optimising the system by minimising cost and maximising EROI.
To better match the methodology presented in Section~\ref{sec:problem-statement} where functions are minimised, $E_{in}$ (i.e. the energy invested in the system) will be used as objective instead of EROI.
The following sections define precisely the two objectives used in the case study.

\subsubsection{System cost}
The first objective is the total annual cost of the system, $f_1 = C_{tot}$, defined as:
\begin{align}
    C_{tot} = \sum_{j\in TECH} \left(\tau(j)C_{inv}(j) + C_{maint}(j)\right) + \sum_{i\in RES} C_{op}(i)\quad.
\end{align}
The yearly system cost is the sum of  $\tau(j)C_{inv}(j)$, the annualised investment cost of each technology with $C_{inv}$ the total investment cost and $\tau$ the annualisation factor, $C_{maint}(j)$, the operating and maintenance cost of each technology and $C_{op}(i)$, the operating cost of the resources.
This last variable is equal to 
\begin{align} \label{equ:c-op}
    C_{op}(i) = \sum_{\substack{t\in T|\{h, td\}\in THTD(t)}} c_{op}(i)\mathbf{F_t}(i, h, td)\quad,
\end{align}
where $c_{op}(i)$ is the cost of resource $i$ in [€/MWh$_\text{fuel}$] and $\mathbf{F_t}(i, h, td)$ corresponds to the use in [MWh$_\text{fuel}$] of resource $i$ at time $(h, td)$.
The values of $c_{op}(i)$ for each resource used in the study case are given in Tables~\ref{tab:resources-data} and \ref{tab:resources-data-res}.
The study of \citet{limpens2019energyscope} or the online documentation \citep{readthedocs} provides more detail on this indicator.\footnote{In the mathematical formulation of the model, an additional factor $t_{op}(h, td)$ is added to equation \eqref{equ:c-op}, \eqref{equ:e-op}, and \eqref{equ:nec-cond-res}.
This parameter is set to 1 in the implementation of the model used in this case study.
It is thus removed from equations for clarity.}

\subsubsection{Energy invested in the system}
The second objective $f_2$ is $E_{in}$, the energy invested in the system over one year:
\begin{align}
    E_{in} = \sum_{j\in TECH} \frac{E_{constr}(j)}{lifetime(j)} + \sum_{i\in RES} E_{op}(i)\quad,
\end{align}
with $E_{constr}(j)$, the energy invested to built technology $j$, annualised by dividing it by its lifetime, and $E_{op}(i)$ the energy to operate, i.e. produce, and transport resource $i$ over one year.
Similarly to the cost indicator, this last variable is equal to 
\begin{align} \label{equ:e-op}
    E_{op}(i) = \sum_{\substack{t\in T|\{h, td\}\in THTD(t)}} e_{op}(i)\mathbf{F_t}(i, h, td)\quad,
\end{align}
where $e_{op}(i)$ is the energy invested (in [MWh/MWh$_\text{fuel}$]) to obtain 1 MWh$_\text{fuel}$ of the resource $i$.
The values of $e_{op}(i)$ for each resource used in the study case are given in Tables~\ref{tab:resources-data} and \ref{tab:resources-data-res}.
More detail is provided by \citet{dumas2022energy} (in which $E_{in}$ is referred to as $E_{in, tot}$).

Minimising $E_{in}$ would be equivalent to maximising EROI, i.e. $E_{out}/E_{in}$, if $E_{out}$, which in our case is the FEC, was constant.
It is not the case in ESTD.
In this model, only the values for the EUD, presented in Table~\ref{tab:eud}, are fixed.
While EUD measures an energy service, FEC measures the quantity of energy used to deliver this service.
FEC is thus always measured in [TWh], while the unit for EUD will depend on the demand.
For instance, the EUD for heat will be measured in [TWh] while [Mt-km] will be used for mobility.
Using technology-dependent conversion factors, FEC can be converted into EUD and vice-versa.
For instance, in ESTD, a FEC of 1 kWh of electricity supplies an EUD of 5.8 passenger-km with a battery-electric car. 
As the conversion factors depend on the installed technologies, which depend on the optimisation results, FEC is an output of the ESTD model and is not constant.
Nonetheless, the constant EUD cannot be employed directly as $E_{out}$ to compute the EROI, as it is an energy service, not an amount of energy.
Therefore, the FEC is used to compute $E_{out}$ and, incidentally, the EROI of the system.

\subsection{Pareto front} \label{sec:case-study-pareto-front}

Once all the elements of the initial optimisation problem \eqref{equ:multi-criteria-initial-problem} are set up, one can compute efficient solutions from the Pareto front using one of the methods described in Section~\ref{sec:appr-pareto-front}.
This case study uses a modified version of the $\epsilon$-constraint method.
It is applied by minimising $E_{in}$ over the feasible space with the additional constraint $C_{tot} \le \epsilon (1 + C^\star_{tot})$ where $\epsilon \in \mathbb{R}_+$ and $C^\star_{tot}$ is the cost-optimal value, i.e. solving
\begin{equation} \label{equ:results-method-epsilon-constraint}
    \begin{split}
        \min_{x \in \mathcal{X} }& \hspace{0.1cm} E_{in} \\
        \text{s.t. } & C_{tot} \le (1+\epsilon)C^\star_{tot}
    \end{split}
\end{equation}
This method is a slight modification of the method described in equation \eqref{equ:method-epsilon-constraint} where $\epsilon$ is a relative rather than absolute value.
It has the benefit of defining the constraint proportionally to the optimal value in the associated objective and thus be directly interpretable.
For instance, if the optimal cost is 75 B€, one would use $\epsilon$ values of 1, 5, and 10\% instead of absolute values of 75.75, 78.75 and 82.5 B€.
To obtain several points over the Pareto front, the method was repeated for different values of $\epsilon$ in $]0, C^e_{tot}/C^\star_{tot}[$ where $C^e_{tot}$ is the value of $C_{tot}$ at the $E_{in}$ optimum and $C^\star_{tot}$ is the cost optimum.

\subsection{Near-optimal spaces} \label{sec:case-study-near-optimal-spacess}
The efficient solutions are used to define approximate near-optimal spaces $\mathcal{X}^{\boldsymbol{\epsilon}}_m$, with $\boldsymbol{\epsilon} = (\epsilon_{C_{tot}}, \epsilon_{E_{in}})$ following equation \eqref{equ:epsilon-space-multi-criteria-def2-approximate-2} of definition \ref{def:epsilon-space-multi-criteria-def2-approximate}.
They are unions of spaces defined around unique, efficient solutions, $\hat{x} \in \mathcal{P}_{\mathcal{X}, m}$.
Each space can be easily defined by adding to the original ESTD model the two linear constraints, which are:
\begin{align}
    C_{tot}(x) \leq (1+\epsilon_{C_{tot}}) C_{tot}(\hat{x})\quad,\\
    E_{in}(x) \leq (1+\epsilon_{E_{in}}) E_{in}(\hat{x})\quad.
\end{align}

\begin{table}
    \renewcommand{\arraystretch}{1.25}
    \centering
    \begin{tabular}{l|r|r} 
        \hline
         & $c_{op}$ & $e_{op}$\\
         & \footnotesize{[€/MWh$_\text{fuel}$]} & \footnotesize{[MWh/MWh$_\text{fuel}$]}\\
        \hline
        \multicolumn{3}{l}{Endogenous resources}\\
        \hline
        Hydro & 0 & 0 \\
        Solar & 0 & 0 \\ 
        Waste & 23.1 & 0.0577 \\ 
        Wet biomass & 5.76 & 0.0559 \\
        Wind & 0 & 0 \\
        Wood & 32.8 & 0.0491\\
        \hline
        \multicolumn{3}{l}{Exogenous resources}\\
        \hline
        Ammonia & 76.0 & $^*0.174$\\ 
        Ammonia (Re.) & 81.8 & $^*0.295$\\
        Diesel & 79.7 & 0.210\\ 
        Bio-diesel & 120 & $^*0.101$ \\
        Elec. import & 84.3 & 0.123 \\
        Gas & 44.3 & 0.0608\\
        Gas (Re.) & 118 & $^*0.147$\\
        Gasoline & 82.4 & 0.281\\ 
        Bio-ethanol & 111 & $^*0.101$\\ 
        H2 & 87.5 & 0.083\\ 
        H2 (Re.) & 119 & $^*0.134$\\
        LFO & 60.2 & 0.204 \\
        Methanol & 82.0 & 0.0798\\ 
        Methanol (Re.) & 111 & $^*0.146$\\ 
    \end{tabular}
    \caption{2035 values for each resource of the input parameters: $c_{op}$, cost of the resource  [€/MWh$_\text{fuel}$] and $e_{op}$, energy invested in obtaining 1 MWh of the resource [MWh/MWh$_\text{fuel}$].
    Most values for $c_{op}$ come from \citep{readthedocs}.
    Data for $e_{op}$ relies on the work of \citet{muyldermansmulti}, who used the \textit{ecoinvent} database \citep{wernet2016ecoinvent}.
    The values preceded by a `$^*$' were updated based on the work by \citet{orban2022energy}.
    Abbreviations: Renewable (Re.), Electricity (Elec.)
    }
    \label{tab:resources-data}
\end{table}
\begin{table}
    \renewcommand{\arraystretch}{1.25}
    \centering
    \begin{tabular}{l|r|r} 
        \hline
        & $c^*_{op}$ & $e^*_{op}$  \\
        & \footnotesize{[€/MWh$_\text{fuel}$]} & \footnotesize{[MWh/MWh$_\text{fuel}$]}\\
        \hline
        Hydro & 53.7 & 0.0489 \\
        Solar & 50.0 & 0.147 \\ 
        Wind & 47.0 & 0.0350\\
    \end{tabular}
    \caption{
    Estimated cost $c^*_{op}$ [€/MWh$_\text{fuel}$] and estimated energy invested in obtaining 1 MWh of the resource $e^*_{op}$ [MWh/MWh$_\text{fuel}$] for three resources: hydro, solar and wind.
    The estimation is done by computing the total cost at the $C_{tot}$ optimum and invested energy at the $E_{in}$ optimum of the technologies that use these resources (i.e. PV for solar, onshore and offshore wind for wind and hydro river for hydro) and then dividing it by the total energy used from these resources at the corresponding optimums, indicated in Table~\ref{tab:optimums-primary-energy}.
    }
    \label{tab:resources-data-res}
\end{table}

\subsection{Necessary conditions} \label{sec:case-study-necessary-conditions}
The last concept to define is the type of necessary conditions computed in the case study.
We are interested in the necessary resources for a transition with a sufficiently good cost and EROI.
We will thus compute the necessary conditions corresponding to the minimum amount of energy that needs to come from a specific individual or group of resources.
Mathematically, the set of such conditions would be:
\begin{align} \label{equ:nec-cond-res}
    \Phi_{\overline{RES}} = \bigg\{\sum_{\substack{i\in \overline{RES},\\ t\in T|\{h, td\}\in THTD(t)}} \mathbf{F_t}(i, h, td) \geq c \bigg\}\quad,
\end{align}
where $\overline{RES} \subseteq RES$ is a set of resources, $F_{i, t}$ the use of resource $r$ at time $t$ and $c \in \mathbb{R}_{+}$.
$\overline{RES}$ can contain any resource.
However, in the context presented in Section~\ref{sec:case-study-context}, we have highlighted a particular interest in two groups of resources: endogenous and exogenous.
We will focus primarily on those two sets and give a more detailed description of their resources.
In ESTD, endogenous resources (noted $RES_{endo}$) include wood, wet biomass, waste, wind, solar, hydro, and geothermal energy. 
Exogenous resources (noted $RES_{exo}$) are the other resources in the model: ammonia, renewable ammonia, imported electricity, methanol, renewable methanol, hydrogen, renewable hydrogen, coal, gas, renewable gas, liquid fuel oil, gasoline, diesel, bio-diesel, and bioethanol. 
Renewable fuels such as renewable ammonia, methanol, and gas are assumed to be produced from renewable electricity.
Tables~\ref{tab:resources-data} and \ref{tab:resources-data-res} list the model's resources and the associated input parameters required to compute the cost and invested energy when employing them.

\section{Results} \label{sec:results}
In this section, we provide the answer to the question that was asked at the beginning of Section~\ref{sec:case-study}:
\begin{center}
    \fbox{\begin{minipage}{0.9\linewidth}
        Which endogenous or exogenous resources are necessary in Belgium to ensure a transition associated with sufficiently good cost and EROI?
    \end{minipage}}
\end{center}
This answer is obtained by computing necessary conditions corresponding to the minimum amount of energy coming from specific resources required to ensure a $\boldsymbol{\epsilon}$-optimality in $C_{tot}$ and $E_{in}$.
However, before diving into the necessary conditions, we first analyse how the system is configured at the two optimums and show the differences between those configurations.
Then, by analysing efficient solutions, we determine how this system evolves when trade-offs are made between $C_{tot}$ and $E_{in}$.
Finally, knowing the Pareto front, we compute $\boldsymbol{\epsilon}$-optimal spaces and examine necessary conditions corresponding to the minimum amount of energy coming from different resources in Belgium.

\subsection{Analysis of the system configuration at the two optimums} \label{sec:results-optimums}
The Belgian energy system is analysed when optimising $C_{tot}$ and $E_{in}$ individually, with a maximum carbon budget $GWP_{tot}$ of 35 Mt\coo-eq/y.
To set a baseline to which we can compare the necessary conditions computed in the following sections, we analyse the amount of endogenous and exogenous resources used at each optimum.
Table~\ref{tab:optimums} shows the value of the two objective functions at the two optimums, and Table~\ref{tab:optimums-primary-energy} details which energy sources are used in the system.

\begin{table}[!ht]
    \centering
    \begin{tabular}{r|r|r|}
        &  $C_{tot}$ & $E_{in}$\\
        & \small{optimum} & \small{optimum}\\
        \hline
        $C_{tot}$ [B\EUR{}/y] & \textbf{52.8} & 56.8\\
        $E_{in}$ [TWh/y] & 74.0 & \textbf{61.0}\\
    \end{tabular}
\caption{Values of $C_{tot}$ and $E_{in}$ objectives at the optimums.}
\label{tab:optimums}
\end{table}
\begin{table}[!ht]
    \begin{tabular}{l|r|r|r}
        Energy &  $C_{tot}$ & $E_{in}$ & Max.\\
        $[\text{TWh/y}]$ & \small{optimum} & \small{optimum} & \small{potential}\\
        \hline
        Endogenous & \textbf{185} & \textbf{164} & \textbf{185}\\
        \hline
        Hydro & 0.469 & 0.486 & $^*0.488$ \\
        Solar & 61.5 & 54.2 & $^*61.6$\\
        Waste & 17.8 & 4.12 & 17.8\\
        Wet Biomass & 38.9 & 38.9 & 38.9\\
        Wind & 42.6 & 43.0 & $^*43.0$\\
        Wood & 23.4 & 23.4 & 23.4\\
        \hline
        Exogenous & \textbf{202} & \textbf{211} & /\\
        \hline
        Ammonia (Re.) & 65.6 & 0 & /\\
        Bio-diesel & 0 & 3.14 & /\\
        Elec. import & 27.6 & 27.6 & 27.6\\
        Gas & 28.2 & 34.5 & /\\
        Gas (Re.) & 4.98 & 48.5 & /\\
        H2 (Re.) & 19.4 & 44.8 & /\\
        Methanol (Re.) & 56.4 & 52.8 & /\\
        \hline
        Total & \textbf{387} & \textbf{375} & /        
    \end{tabular}
    \caption{Amount of energy used from each endogenous and exogenous resource at $C_{tot}$ and $E_{in}$ optimums.
    The last column shows the maximum potential of each resource.
    The potentials preceded by a `*' are computed from the maximum capacity and capacity factors of the technologies using these resources.
    The other potentials are directly fixed as parameters.
    The resources for which a `/' is noted have unlimited potential.
    }
    \label{tab:optimums-primary-energy}
\end{table}

\subsubsection{Results at the $C_{tot}$ optimum}
The optimal cost $C^\star_{tot}$ is equal to 52.8 B\EUR{}/y.
At this optimum, the total amount of primary energy used in the system is 387 TWh/y, 48\% of which comes from endogenous resources and the rest from exogenous resources.

For endogenous resources, the values for wet biomass, waste and wood are equal to their maximum potentials - set as input model parameters.
This observation makes sense as the $c_{op}$ of these resources in Table~\ref{tab:resources-data} indicate they are among the cheapest resources.
The hydro, solar and wind energy quantities are also very close to their maximum potential.
For these resources, the maximum is not set directly on the quantity of energy but on the capacities of the technologies using these resources.
For instance, the model can install a maximum of 6 GW of offshore wind turbines and 10 GW of onshore wind turbines, which are the two technologies using wind as a resource.
These maximum capacities can then be multiplied by the capacity factors of the corresponding technologies to obtain a maximum energy potential.
Moreover, these resources are considered free in terms of cost and invested energy, as shown in Table~\ref{tab:resources-data}.
The cost of using them arises from the technologies to extract them from the environment.
Table~\ref{tab:resources-data-res} shows approximated values for $c_{op}$ and $e_{op}$.
They are computed by dividing the cost or energy invested for building and maintaining the technologies using them by the total energy used from these resources - shown in Table~\ref{tab:optimums-primary-energy}.
These approximated values show that hydro, solar and wind are among the cheapest resources, which explains their extensive use.

The model has no maximum potential for exogenous resources except for imported electricity.
This potential is reached as, even though $c_{op}$ is relatively high for imported electricity, it does not require any conversion technology to produce the final electricity demand.
Some 65.6 TWh/y of renewable ammonia is used in the system, 55.4 TWh/y of which is used for electricity production and low-temperature heat generation, while the remaining 10.2 TWh/y is used to satisfy non-energy demand.
Most renewable methanol is used to produce high-value chemicals, even though 3.6 TWh/y of this resource is used for fuelling boat freight.
Finally, gas (renewable or not) is used to produce heat and electricity and fuel buses for public mobility.

\subsubsection{Results at the $E_{in}$ optimum}
The optimal energy invested  $E^\star_{in}$ amounts to 61 TWh/y.
Among the 375 TWh/y of primary energy in the system, 164 TWh/y (44\%) come from endogenous resources and 211 TWh/y (56\%) from exogenous resources.
A series of resources, including wet biomass, wood, wind, hydro and imported electricity, are used at or near their maximum potential.
This is not the case for waste and solar.
In particular, for solar, $e_{op}$ is about three times higher than any other endogenous resource.
This result can be explained by the higher energy needed to build 1 GW of PV combined with a low average capacity factor compared to hydro river plants or wind turbines.
Some electricity is produced using natural and renewable gas, while ammonia for non-energy demand is produced from H2 using the Haber-Bosch process.
The remaining amount of gas is used to produce heat.
Finally, high-value chemicals are produced using renewable methanol, while 3.14 TWh/y of bio-diesel is used for boat freight.

\subsubsection{Comparison}
Table~\ref{tab:optimums} shows how the two objective functions vary from one optimum to the other.
The increase in cost when optimising $E_{in}$ is limited to 7.67\%.
Invested energy at the $C_{tot}$ optimum is around 74 TWh/y, representing an increase of more than 20\% from the $E_{in}$ optimal value.

As shown in Table~\ref{tab:optimums-primary-energy}, the total amount of energy needed in the system differs only by 3\%, but there are some differences between the two energy mixes.
At the $C_{tot}$ optimum, the energy coming from endogenous resources is 21 TWh/y higher, while energy from exogenous resources is 9 TWh/y smaller.
At each optimum, the share of endogenous resources in the energy mix (48\% and 44\%, respectively) is close to the maximum of 42\% primary energy coming from endogenous resources computed by \citet{limpens2020belgian}. 
These values confirm the substantial dependence of Belgium on imported resources to supply its energy consumption.

Looking at individual resources, solar and renewable ammonia, used to produce electricity when optimising $C_{tot}$, are replaced by fossil and renewable gas at the $E_{in}$ optimum.
At this optimum, a percentage of the total 80 TWh/y of gas is used to produce high-temperature heat instead of waste.
The additional 35.4 TWh/y of renewable hydrogen is used for three things: ammonia production (which was directly imported when optimising cost), combined heat and electricity production, and public mobility.
Finally, while boat freight is fuelled using renewable methanol at the $C_{tot}$ optimum, bio-diesel is the preferred option at the $E_{in}$ optimum.

\begin{figure}
    \centering
    \includegraphics[width=\linewidth]{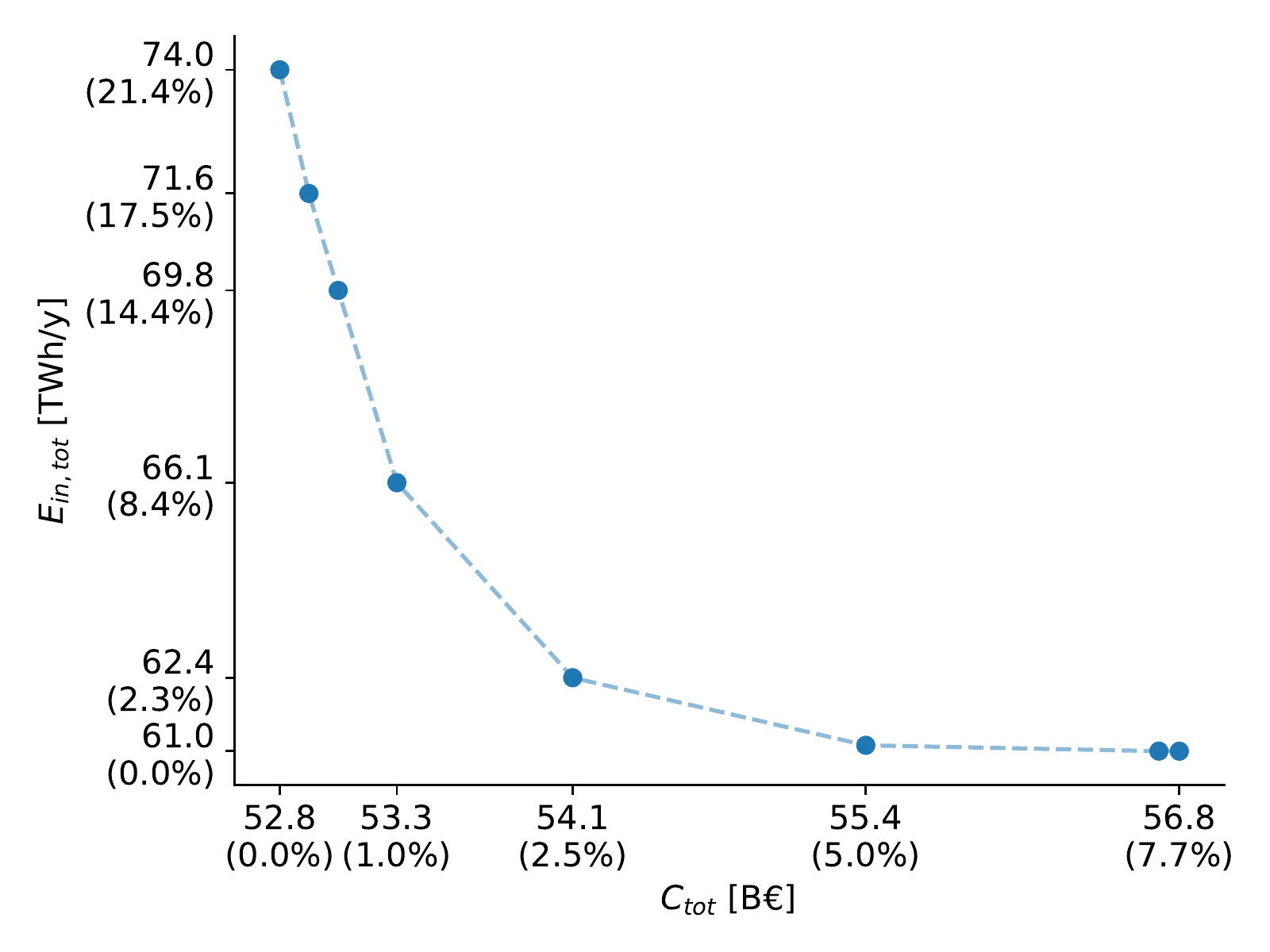}%
    \caption{Approximated Pareto front showing trade-offs between $C_{tot}$ and $E_{in}$.
    On the axis, the absolute values of $C_{tot}$ and $E_{in}$ are shown and completed, in parenthesis, by the deviations from the optimal objective values in each objective.
    For instance, for $C_{tot}$, the value $C_{tot}/C^\star_{tot}-1$ is shown in parenthesis.
    }
    \label{fig:pareto_front}
\end{figure}

\subsection{Pareto Front}

\begin{figure*}
    \centering
    \subfloat[Endogenous resources
    ]
    {\includegraphics[width=3.3in]{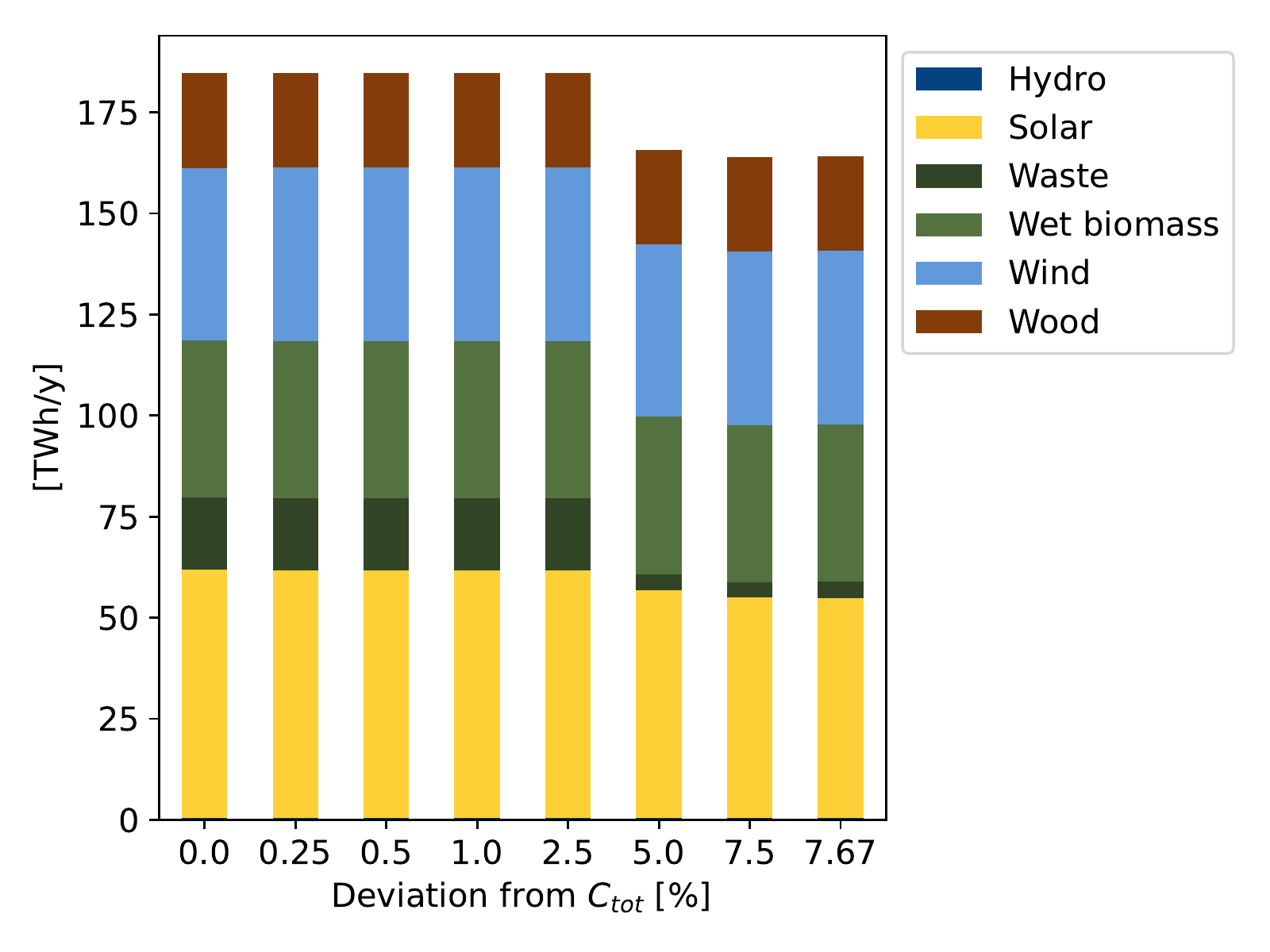}
    \label{fig:pareto_locals_bar}}
    \subfloat[Exogenous resources]
    {\includegraphics[width=3.3in]{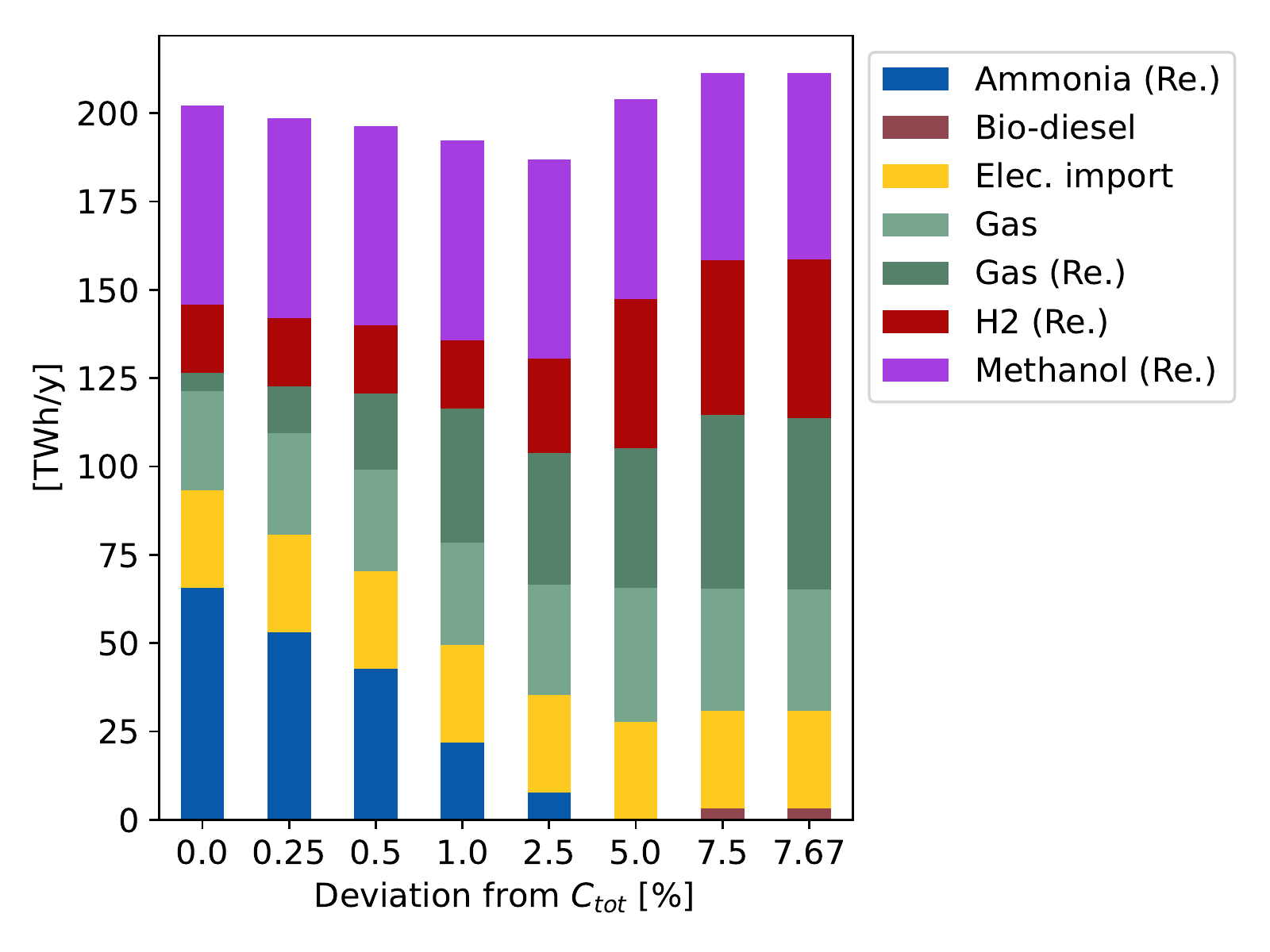}%
    \label{fig:pareto_imports_bar}}
    \caption{Energy [TWh/y] coming from (a) endogenous and (b) exogenous resources at efficient solutions representing different trade-offs between $C_{tot}$ and $E_{in}$. 
    The leftmost bars show these values at the $C_{tot}$ optimum, while the rightmost bar shows these values at the $E_{in}$ optimum.
    The bars in the middle are characterized by their deviation in [\%] from the $C_{tot}$ optimum.
    Abbreviations: Renewable (Re.), Electricity (Elec.).
    }
    \label{fig:pareto_bars}
\end{figure*}

Figure~\ref{fig:pareto_front} shows the values of $C_{tot}$ and $E_{in}$ at the efficient solutions obtained using the method described in Section~\ref{sec:case-study-pareto-front} for values of $\epsilon$ equal to 0.25, 0.5, 1.0, 2.5, 5.0 and 7.5\%.
The two additional points at the curve extremes correspond to each objective's optimum.
The axes are labelled both in terms of the absolute values of the objective functions but also - in parenthesis - in terms of the deviations of these values from the optimal objective value, i.e. $C_{tot}/C^\star_{tot}-1$ and $E_{in}/E^\star_{in}-1$.

This graph shows that $E_{in}$ decreases quite rapidly, saving 10 TWh/y out of 74 TWh/y ($\sim$ -14\%) when increasing $C_{tot}$ by a relatively small amount of 2.5\%.
This behaviour can also be interpreted as: choosing the optimal cost implies a considerable addition in invested energy.
Inversely, as already mentioned, $C_{tot}$ is still relatively low at the $E_{in}$ optimum, i.e. it only increases by 7.5\%.

Figure~\ref{fig:pareto_bars} shows the amount of endogenous and exogenous resources used at each efficient solution, starting on the left with the cost optimum and moving towards the invested energy optimum on the right.
As stated when comparing optimums, there is only a minor change for endogenous resources when going from one optimum to the other.
This change, the reduction of solar and waste energy, appears when allowing a 5\% deviation in cost.

More change is happening for exogenous resources (Figure~\ref{fig:pareto_imports_bar}).
As we increase cost and decrease the invested energy, ammonia is gradually replaced by gas (both natural and renewable).
At a 2.5\% cost increase, the amount of renewable H2 starts increasing.
Ammonia is wholly removed from the system at 5\%, while natural gas use reaches its maximum and starts to decline.
The same happens for renewable ammonia when reaching a 7.5\% cost increase, and some bio-diesel appears.
Overall, the change in the total amount of exogenous resources used is non-monotonic.
Starting to decrease, it then increases when reaching the 5\% threshold, corresponding to the drop in endogenous resources use.

\subsection{Necessary conditions} \label{sec:results-necessary-conditions}

\begin{figure*}[!ht]
    \centering
    \subfloat[Endogenous resources
    ]
    {\includegraphics[width=3in]{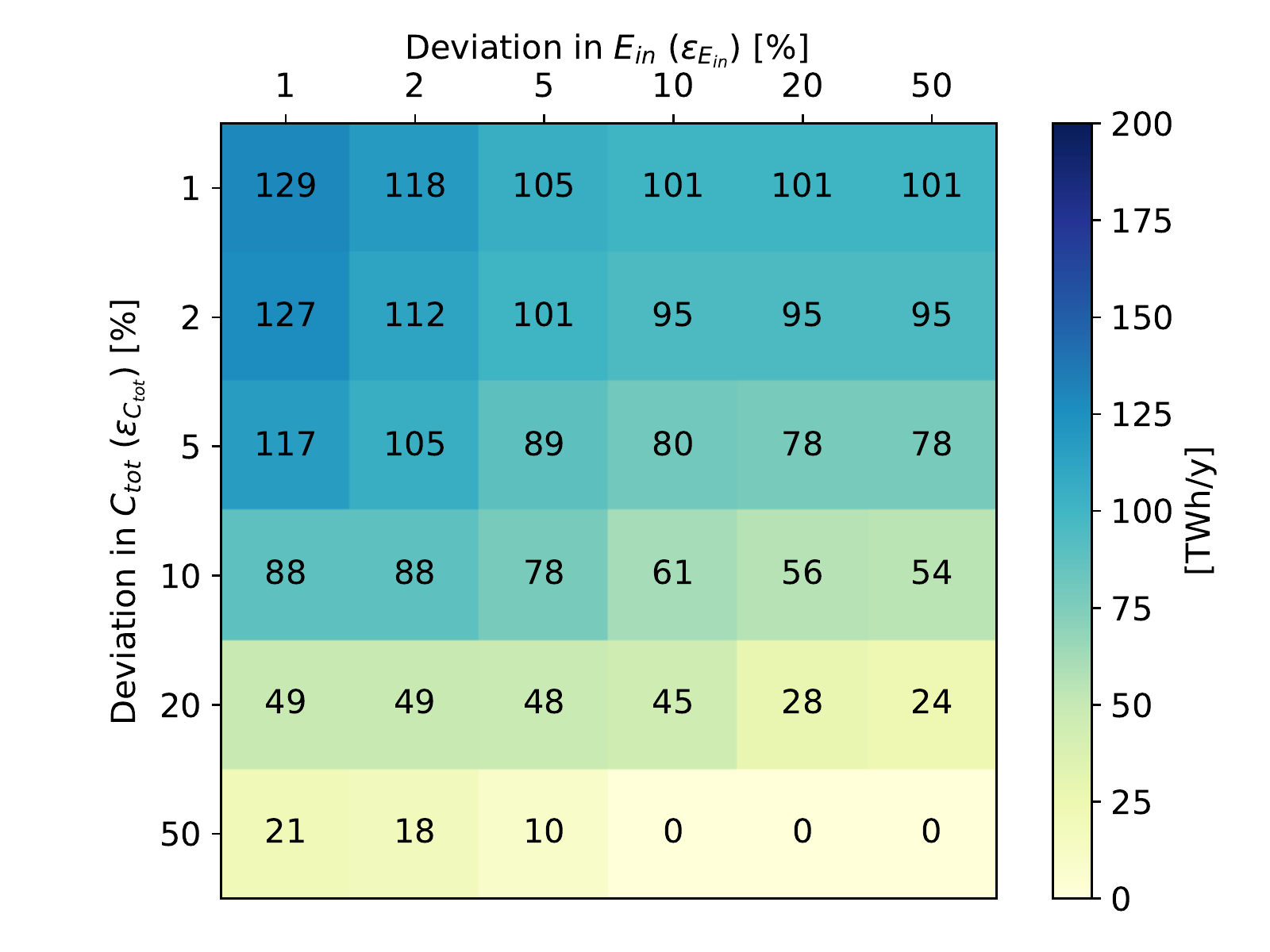}
    \label{fig:necessary_conditions_locals}}
    \subfloat[Exogenous resources]
    {\includegraphics[width=3in]{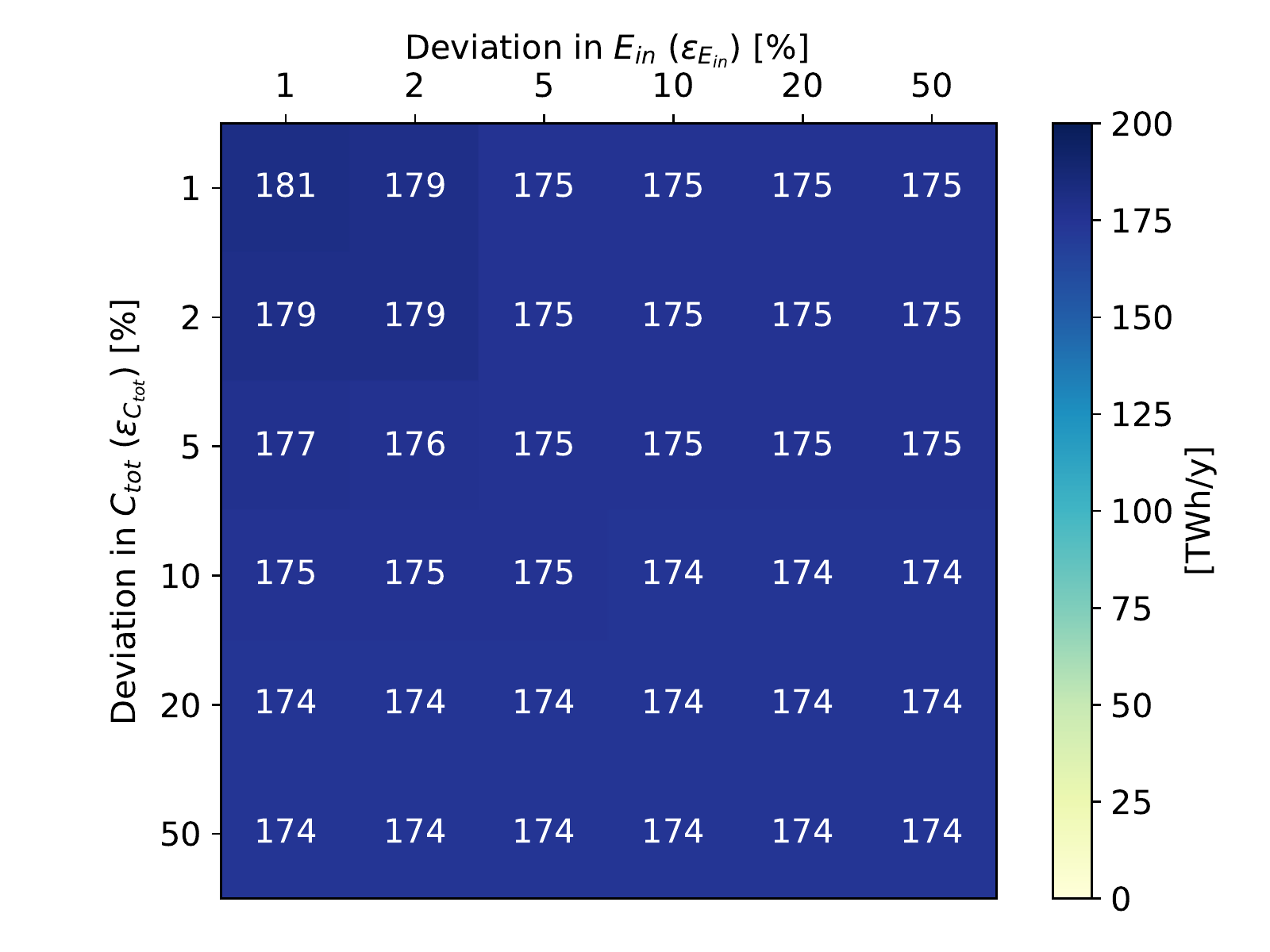}%
    \label{fig:necessary_conditions_imports}}\\
    \subfloat[Renewable methanol]
    {\includegraphics[width=\textwidth/3]{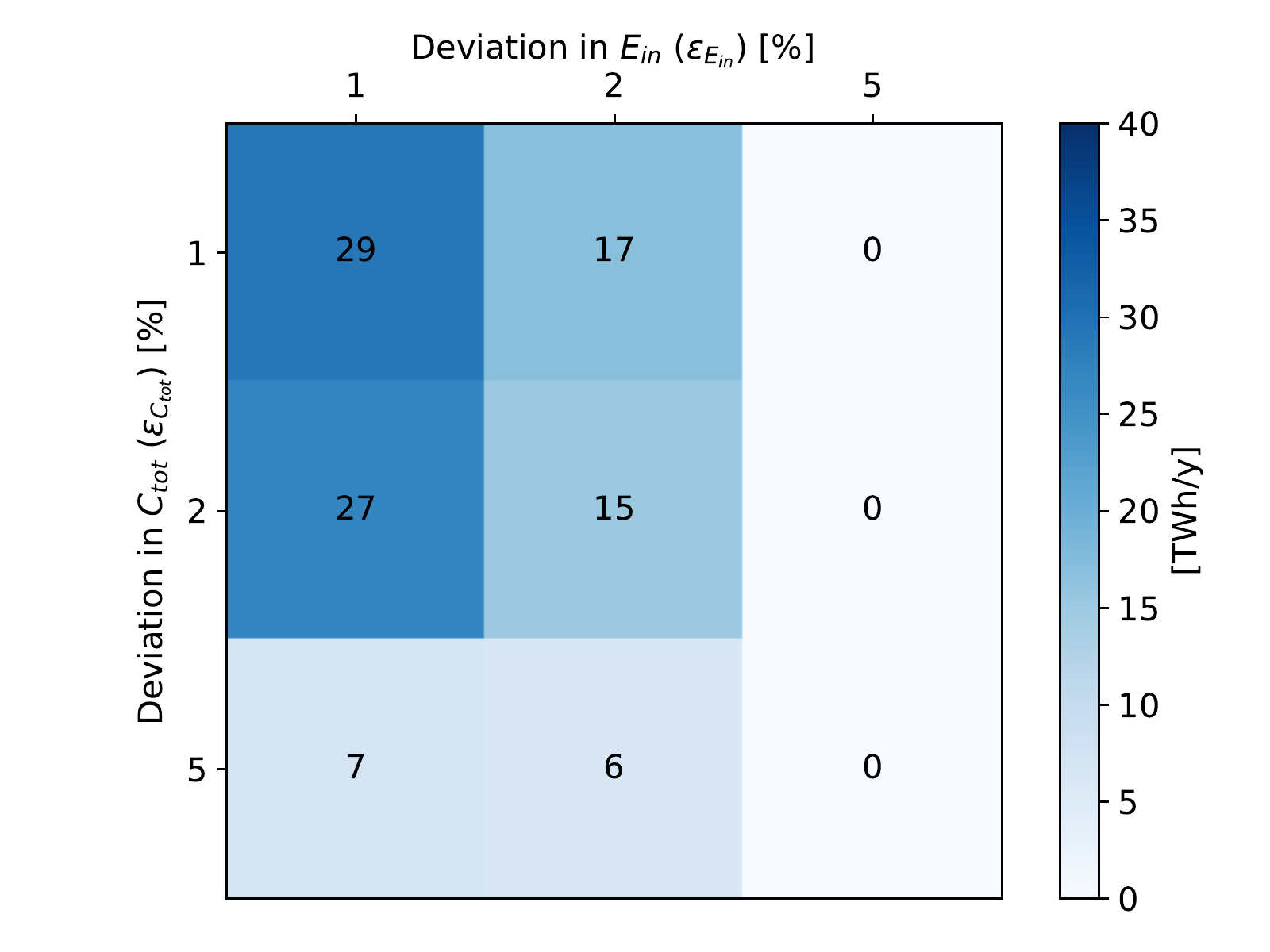}
    \label{fig:necessary_conditions_methanol_re}}
    \subfloat[Gas]
    {\includegraphics[width=\textwidth/3]{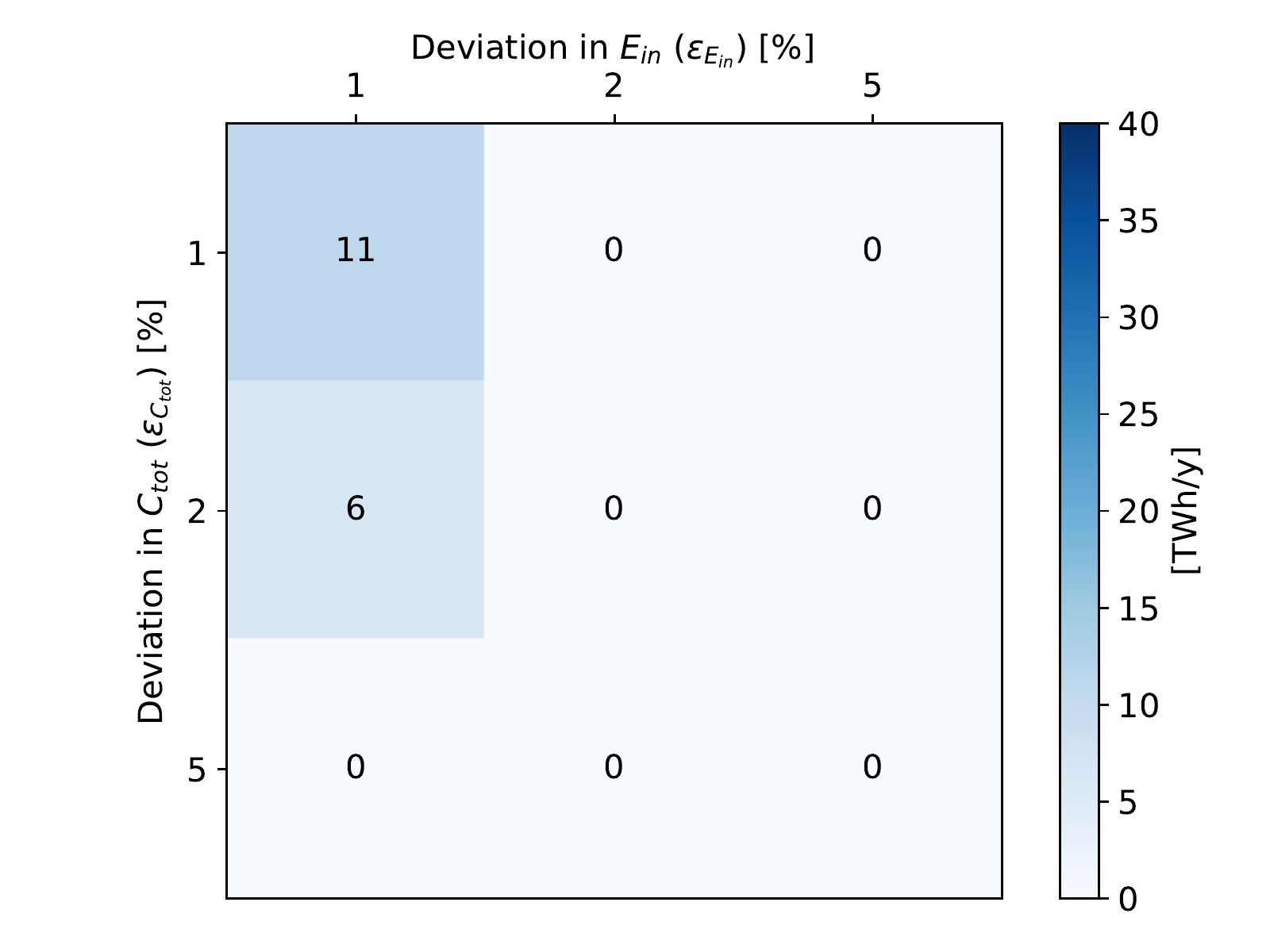}%
    \label{fig:necessary_conditions_gas}}
    \subfloat[Electricity import]
    {\includegraphics[width=\textwidth/3]{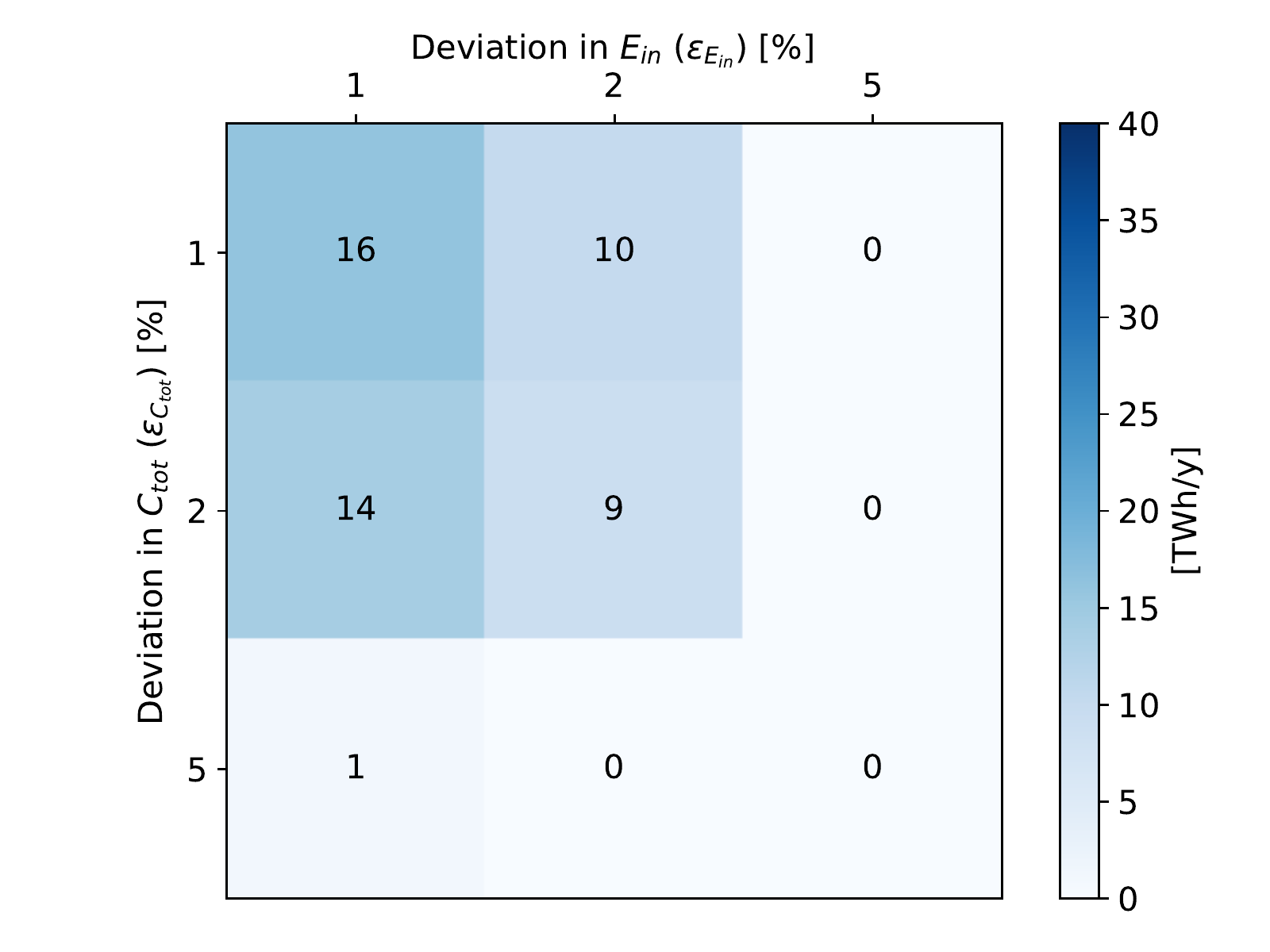}%
    \label{fig:necessary_conditions_elec_import}}
    \caption{Values $\Tilde{c}$ of necessary conditions (in [TWh/y]) for conditions of type $\Phi_{\overline{RES}}$.
    The set of resources $\overline{RES}$ corresponds to endogenous resources $RES_{endo}$ and exogenous resources $RES_{exo}$ for graph (a) and (b), respectively, while for graphs (c), (d) and (e), this set resumes to a single resource: renewable methanol, gas and imported electricity, respectively.
    The values correspond to the minimum amount of energy that needs to come from these sets of resources to ensure a constrained deviation in $C_{tot}$ and $E_{in}$.
    These deviations are defined by the suboptimality coefficients vector $\boldsymbol{\epsilon} = (\epsilon_{C_{tot}}, \epsilon_{E_{in}})$.
    For (a) and (b), all the combinations of the following percentages are taken as coefficients vectors: 1\%, 2\%, 5\%, 10\%, 20\%, and 50\%.
    For (c), (d) and (e), they are limited to the combinations of 1\%, 2\%, and 5\%.
    }
    \label{fig:necessary_conditions}
\end{figure*}

Analysing efficient solutions gives a first appreciation of the variety of system configurations, offering a trade-off between different objectives.
However, using the necessary conditions, we can go one step further by providing features respected by all those solutions and some slightly less efficient solutions.
We use Algorithm~\ref{alg:cap-2} to compute non-implied necessary conditions stemming from different sets of conditions of the type defined by \eqref{equ:nec-cond-res} in Section~\ref{sec:case-study-necessary-conditions}.
The main parameter defining these conditions is $\overline{RES}$, the set of resources over which the constrained sum is computed.
The output of this algorithm is a value $\Tilde{c}$, which defines a non-implied necessary condition for this set of resources.
Practically, this value represents the minimum amount of energy that needs to come from this set of resources to ensure that $C_{tot}$ and $E_{in}$ do not deviate by more than an $\boldsymbol{\epsilon}$ fraction from at least one solution in the Pareto front.
We will first compute this $\Tilde{c}$ value for conditions defined using the set of endogenous and the set of exogenous resources.
We will then look at sets containing one individual resource.

\subsubsection{Endogenous vs exogenous resources}
In this first section, we compare the values $\Tilde{c}$ of non-implied necessary conditions computed from the sets $\Phi_{RES_{endo}}$ and $\Phi_{RES_{exo}}$.
These conditions are computed for different values of deviations $\boldsymbol{\epsilon}$.
In this case, the tuples $\boldsymbol{\epsilon} = (\epsilon_{C_{tot}}, \epsilon_{E_{in}})$ corresponds to all the possible combinations of 1\%, 2\%, 5\%, 10\%, 20\%, and 50\%.

Comparing Figures~\ref{fig:necessary_conditions_locals} and \ref{fig:necessary_conditions_imports} shows how the behaviours of the minima in endogenous and exogenous resources are very different.
For endogenous resources, the minimum for deviations of 1\% in both objectives is already down to 130 TWh/y, representing a 42\% and 26\% decrease from the $C_{tot}$ and $E_{in}$ optima, respectively.
This amount is divided by more than two when the deviation reaches 10\% in both objectives, leaving only 60 TWh/y left from endogenous resources.
The $\Tilde{c}$ value then reaches 0 TWh/y when allowing for an increase of 50\% in $C_{tot}$.
These results show that energy from endogenous resources can be reduced by a significant amount for reasonably low increases in cost and invested energy.

For exogenous resources, there is little to no decrease in the total energy needed.
Starting from 202 and 211 TWh/y at the optimums in cost and energy invested, the minimum amount of this type of energy is still around 174 TWh/y (i.e. -20\% and -15\% respectively) for deviations of 10\%.
Most of the decrease is already present for deviations of 1\% with an amount of energy of 180 TWh/y, which is only 6 TWh/y less than the energy used at one of the efficient solutions.
The $\Tilde{c}$ value of non-implied necessary conditions then plateaus at 174 TWh/y.
This result shows how, contrarily to endogenous resources, exogenous resources are essential, whatever the cost and energy invested.
Indeed, to respect a $GWP_{tot}$ constraint of 35 Mt\coo-eq/y, at least 174 TWh/y of energy needs to be imported.

\subsubsection{Individual exogenous resources}
We have shown that a certain amount of exogenous resources is necessary due to limited endogenous resources.
However, the previous results do not show which specific exogenous resource is essential.
This analysis can be done by computing necessary conditions for groups of conditions $\Phi_{\{i\}}$ where $i \in RES$ corresponds to an individual resource.
We could perform this analysis for all individual resources, but in Figure~\ref{fig:pareto_imports_bar}, the amounts of renewable methanol, gas, and imported electricity are quasi-constant across the Pareto front.
Therefore, it is interesting to focus on these resources to see if they are essential or if we can eliminate them by increasing the cost or the invested energy.
In this section, we compute non-implied necessary conditions corresponding to the minimum energy from these three resources.

The $\Tilde{c}$ values of non-implied necessary conditions for deviations $\epsilon_{C_{tot}}$ and $\epsilon_{E_{in}}$ of 1\%, 2\% and 5\% are shown in Figures~\ref{fig:necessary_conditions_methanol_re}, \ref{fig:necessary_conditions_gas}, and \ref{fig:necessary_conditions_elec_import}.
We limit the analysis to deviations of 5\% as we can see that we are already equal (or near to) 0 TWh/y for all three resources at this percentage.
The amount of energy coming from the resources at the $C_{tot}$ and $E_{in}$ optimums are respectively 56.4 and 52.8 TWh/y for renewable methanol, 28.2 and 34.5 for gas, and 27.6 (at both optima) for imported electricity.
The minimum energy from each resource is around 50\% lower than at the efficient solutions when allowing deviations of 1\% in each objective.
For renewable methanol and gas, the amount of necessary energy is more sensitive to deviations in invested energy than to deviations in cost.
However, the conclusion is similar for the three resources: for a relatively small increase in cost and invested energy, they can be replaced by other resources.

\section{Conclusion} \label{sec:conclusion}

This paper presents an extension of the concepts of necessary conditions for epsilon-optimality introduced by \citet{DUBOIS2022108343} for mono-objective optimisation to multi-objective optimisation.
This methodology provides a new way of overcoming the limitation that energy system models face when focusing on cost.

We first highlight the limits of cost-based optimisation problems and a series of existing solutions: scenario analysis, multi-objective optimisation and near-optimal spaces analysis.
To outline how this last methodology can be combined with multi-objective optimisation, the concepts of epsilon-optimality and necessary conditions are presented in a single-objective setup and extended to multi-objective problems.
Then, these concepts are illustrated in a case study addressing the Belgian energy transition and answering the question: ``Which endogenous or exogenous resources are necessary in Belgium to ensure a transition associated with sufficiently good cost and EROI?".
The answer is obtained by computing non-implied necessary conditions corresponding to the minimum amount of energy coming from different sets of resources to ensure a constrained deviation in cost and energy invested.
The results show Belgium's dependence on imported resources but that no individual resource is essential.

This paper introduces a methodology and applies it for the first time to a concrete problem.
There are, therefore, limitations to overcome and work tracks to explore.
Regarding methodology, the main limitation is that, as explained in Section~\ref{sec:methodology-multi}, in a multi-objective optimisation problem, the epsilon-optimal space can only be approximated.
This approximation, in turn, removes the guarantees of finding non-implied necessary conditions.
Increasing the number of efficient solutions used is a straightforward way to improve the quality of the results, but it comes at the expense of computational time.
More research could be done to evaluate how the number and spread of efficient solutions affect the results.
Taking a more general viewpoint, the techniques presented in this paper provide insights using fixed feasible spaces and objective functions.
To reach their full potential, they ought to be combined with methods for fighting parametric uncertainty, such as sensibility analysis.
We explored how the concept of necessary conditions could be extended to a multi-objective set-up.
However, other methodologies were developed to explore near-optimal spaces in a mono-objective setup as presented, for instance, in \cite{PRICE2017356, li2017investment, pedersen2021modeling, nacken2019integrated}.
An interesting research track would thus be to extend these methodologies in a multi-objective setup.

To illustrate the methodology, a specific case study with a limited scope was used.
A natural extension is thus to replicate this case study for other countries or territories where limited resources and energy dependence represent a challenge.
Alternative questions regarding other resources or technologies could also be addressed using this framework.
Finally, going one step further, a follow-up to this work would involve studying near-optimal spaces for other optimisation problems inside or outside the energy systems field.

\section*{Acknowledgements}

Antoine Dubois is a Research Fellow of the F.R.S.-FNRS, of which he acknowledges the financial support.

\section*{Data Availability}

Datasets related to this article can be found at \url{https://zenodo.org/record/7665440#.Y_YyltLMIUE}, hosted at Zenodo (\cite{dubois2023code}).

\bibliographystyle{elsarticle-num-names} 
\bibliography{main}

\appendix
\newpage
\section{Example} \label{sec:appendix-example}

\noindent The functions depicted in Figures~\ref{fig:epsilon-optimal-space-uni}, \ref{fig:epsilon-optimal-space-2d}, and \ref{fig:epsilon-space-approximate-pareto} are 
\begin{align*}
    &f_1(x) = 10*(2x-0.75)^2+2\quad\text{and}\\\
    &f_2(x) = 10*(x-0.75)^2+1.5\quad.
\end{align*}
The coordinates of their minimums are $(x^*_1, f_1(x^*_1)) = (0.375, 2)$ and $(x^*_2, f_2(x^*_2)) = (0.75, 1.5)$, respectively.

\subsection*{One-dimensional epsilon-optimal space}
In Figure~\ref{fig:epsilon-optimal-space-uni}, the $\epsilon$-optimal space $\mathcal{X}^\epsilon$ of a one-dimensional optimisation problem was obtained by first computing 
\begin{align*}
    (1+\epsilon_1)f_1(x^*_1) = (1+0.25)*2 = 2.5   
\end{align*}
where $\epsilon_1 = 0.25$.
Then, the limits of $\mathcal{X}^\epsilon$ can be obtained by computing the inverse image of this value, i.e. the set $\{0.263, 0.487\}$, which leads to $\mathcal{X}^\epsilon = [0.263, 0.487]$.

\subsection*{Two-dimensional epsilon-optimal space}
In Figure~\ref{fig:epsilon-optimal-space-2d-full} and \ref{fig:epsilon-optimal-space-2d-part}, the Pareto front $\mathcal{P}_\mathcal{X}$ is represented in green.
This set of points respects definition~\ref{def:pareto-front} of a Pareto front.
Indeed, each point $x$ in the interval $[x^\star_1, x^\star_2]$ is such that $\not\exists \hat{x}\in\mathcal{X}$ where $f_1(\hat{x}) < f_1(x)$ and $f_2(\hat{x}) < f_2(x)$.

In Figure \ref{fig:epsilon-optimal-space-2d-part}, a subset of the $\boldsymbol{\epsilon}$-optimal space $\mathcal{X}^{\boldsymbol{\epsilon}}$ of a two-dimensional optimisation problem is computed for a suboptimality coefficients vector $\boldsymbol{\epsilon} = (\epsilon_1, \epsilon_2) = (0.25, 0.6)$.
This subset is computed from the point $\hat{x} = 0.6$, which is part of $\mathcal{P}_\mathcal{X}$.
To obtain the subset of $\mathcal{X}^{\boldsymbol{\epsilon}}$, the images of $\hat{x}$, $f_1(\hat{x}) = 4.025$ and $f_2(\hat{x}) = 1.725$, are computed.
Multiplying these values by the corresponding suboptimality coefficients gives
\begin{align*}
    &(1+\epsilon_1)f_1(\hat{x}) = (1+0.25)*4.025 = 5.03\quad\text{and}\\
    &(1+\epsilon_2)f_2(\hat{x}) = (1+0.6)*1.725 = 2.76\quad.    
\end{align*}

The inverse image of these values are $\{0.0997, 0.65\}$ for $f_1$ and $\{0.395, 1.105\}$ for $f_2$.
The set of points respecting $\forall k\; f_k(x) \leq (1+\epsilon_k) f_k(\hat{x})$ are then contained in $[0.395, 0.65]$.

To obtain the full $\boldsymbol{\epsilon}$-optimal space depicted in Fig.~\ref{fig:epsilon-optimal-space-2d-full}, one should repeat this process with all points in $\mathcal{P}_{\mathcal{X}}$.
However, in this simple example, one can easily compute the limits of the complete space by using the two optimums, which are the extreme points of the Pareto front.
These limits are obtained by taking the inverse images of 
\begin{align*}
    &(1+\epsilon_1)f_1(x^*_1) = (1+0.25)*2 = 2.5\quad\text{and}\\
    &(1+\epsilon_2)f_2(x^*_2) = (1+0.6)*1.5 = 2.4\quad,
\end{align*}
which gives $\{0.263, 0.487\}$ and $\{0.45, 1.05\}$.
The lower and upper bound of $\mathcal{X}^{\boldsymbol{\epsilon}}$ are then respectively given by the lower and upper bound of those two sets, i.e. $\mathcal{X}^{\boldsymbol{\epsilon}} = [0.263, 1.05]$.

\subsection*{Approximate Pareto fronts and epsilon-optimal spaces}
Figures~\ref{fig:epsilon-space-approximate-pareto-3}, \ref{fig:epsilon-space-approximate-pareto-2}, and \ref{fig:epsilon-space-approximate-pareto-1} show approximate $\boldsymbol{\epsilon}$-optimal spaces for three different set of efficient points.
These sets are 
\begin{enumerate}
    \item Fig.~\ref{fig:epsilon-space-approximate-pareto-3}: [(2.0, 2.91), (3.41, 1.85), (7.62, 1.5)];
    \item Fig.~\ref{fig:epsilon-space-approximate-pareto-2}: [(2.9, 2.01), (2.99, 1.97), (3.09, 1.94), (3.19, 1.91), (3.3, 1.88), (3.41, 1.85), (3.52, 1.82), (3.64, 1.8), (3.77, 1.77), (3.9, 1.75), (4.03, 1.72)];
    \item Fig.~\ref{fig:epsilon-space-approximate-pareto-1}: [(2.0, 2.91), (2.06, 2.64), (2.23, 2.40), (2.51, 2.19), (2.9, 2.01), (3.41, 1.85), (4.03, 1.72), (4.76, 1.63), (5.61, 1.56), (6.57, 1.51), (7.62, 1.5)]. 
\end{enumerate}

\section{Types of end-use demand in EnergyScope-TD} \label{sec:appendix-eud}
Four main types of EUD are considered in the model: electricity, heat, transport, and non-energy demand.
Electricity is further divided between lighting and other electricity uses. 
Heat is subdivided into high-temperature heat for industry, low temperature for space heating, and low temperature for hot water. 
Mobility is composed of public and private passenger mobility and freight demands.
Finally, the non-energy demand includes demand for ammonia, methanol, and high-value chemicals (HVCs).
Table \ref{tab:eud} lists the values for each EUD type in 2035 based on \citet{limpens2021generating}.
\begin{table}[!ht]
\centering
\begin{tabular}{l|l|r} 
    \hline
    EUD type & Unit & EUD \\
    \hline
    Electricity (other) & TWhe & 62.1\\
    Lighting & TWhe & 30.0\\
    Heat high T. & TWh & 50.4\\ 
    Heat low T. (SH) & TWh & 118\\ 
    Heat low T. (HW) & TWh & 29.2\\
    Passenger mobility & Mpass.-km & 194 \\
    Freight & Mt-km & 98.0 \\
    Non-energy & TWh & 53.1\\   
\end{tabular}
\caption{2035 Belgian end-use demand (EUD) value by type based on \citet{limpens2021generating}. Abbreviations: temperature (T.), space heating (SP), hot water (HW), passenger (pass.).}
\label{tab:eud}
\end{table}

\end{document}